\documentclass[12pt,namelimits,sumlimits]{amsart}
\usepackage{citesort}
\usepackage{amssymb,amsmath}
\usepackage[mathscr]{eucal}

\usepackage{psfrag}

\usepackage{graphicx}
\usepackage{amsmath,amsthm}
\usepackage{graphics}
\usepackage{color}
\usepackage{epsfig}
\usepackage{fullpage}
\usepackage{amssymb,amsmath}
\usepackage[mathscr]{eucal}

\newcounter{mtheorem}

\newtheorem{theorem}{Theorem}[section]
\newtheorem{lemma}[theorem]{Lemma}
\newtheorem{prop}[theorem]{Proposition}

\newtheorem{notation}[theorem]{Notation}
\newtheorem{convention}[theorem]{Convention}
\newtheorem{question}[theorem]{Question}

\newtheorem{definition}[theorem]{Definition}

\theoremstyle{remark}
\newtheorem{remark}[theorem]{Remark}

\numberwithin{equation}{section}

\newcommand{\tb}{{\underline{\tau}}}

\newcommand{\RRR}{\mathsf{R}}

\newcommand{\skernel}{\mathscr{K}}

\newcommand{\beq}{\begin{equation}}
\newcommand{\eeq}{\end{equation}}
\newcommand{\bea}{\begin{eqnarray}}
\newcommand{\eea}{\end{eqnarray}}

\newcommand{\Rbb}{\mathbb{R}}

\newcommand{\Zbb}{\mathbb{Z}}
\newcommand{\Nbb}{\mathbb{N}}

\newcommand{\Sph}{\mathbb{S}}

\newcommand{\cunder}{\underline{c}\,}

\newcommand{\xibold}{{\boldsymbol{\xi}}}

\DeclareMathOperator*{\avg}{avg}

\newcommand{\Lh}{\ensuremath{\mathcal L_h}}

\newcommand{\Tbb}{\mathbb{T}}

\newcommand{\vecK}{{\vec{K}}}
\newcommand{\dd}{{\boldsymbol{\partial}\mspace{.8mu}\!\!\!\!\boldsymbol{/}\,}}

\newcommand{\tran}{{T\mspace{0.9mu}\!\!\!\!\boldsymbol{/}\,}}

\newcommand{\W}{{\mathcal{W}}}
\newcommand{\N}{{\mathcal{N}}}
\newcommand{\C}{{\mathcal{C}}}
\newcommand{\Cunder}{\underline{{\mathcal{C}}}}
\newcommand{\What}{{\widehat{{\mathcal{W}}}}}
\newcommand{\Chat}{{\widehat{{\mathcal{C}}}}}
\newcommand{\X}{\boldsymbol{X}}
\newcommand{\QQ}{\boldsymbol{Q}}
\newcommand{\A}{\boldsymbol{A}}
\newcommand{\E}{\boldsymbol{E}}
\newcommand{\CC}{\boldsymbol{C}}
\newcommand{\cc}{\boldsymbol{c}}
\newcommand{\ttt}{\boldsymbol{t}}
\newcommand{\Z}{\boldsymbol{Z}}
\newcommand{\WW}{\boldsymbol{W}}
\newcommand{\Zti}{\widetilde{\boldsymbol{Z}}}
\newcommand{\xibo}{{\boldsymbol{\xi}}}
\newcommand{\al}{\underline{\boldsymbol{\alpha}}}
\newcommand{\cch}{\boldsymbol{\chi}}
\newcommand{\llambda}{\boldsymbol{\lambda}}
\newcommand{\cchhat}{\widehat{\boldsymbol{\chi}}}
\newcommand{\alhat}{\widehat{\boldsymbol{\alpha}}}

\newcommand{\beche}{\check{\boldsymbol{\beta}}}
\newcommand{\thetahat}{\widehat{\theta}}
\newcommand{\bb}{\boldsymbol{b}}
\newcommand{\bbe}{\boldsymbol{\beta}}
\newcommand{\bbeha}{\widehat{\boldsymbol{\beta}}}
\newcommand{\bbeti}{\widetilde{\boldsymbol{\beta}}}
\newcommand{\zzero}{{\boldsymbol{0}}}
\newcommand{\BB}{\widetilde{\boldsymbol{b}}}
\newcommand{\ffti}{\widetilde{\boldsymbol{f}}}
\newcommand{\gm}{\boldsymbol{\nu}}
\newcommand{\Xhat}{\widehat{\boldsymbol{X}}}

\newcommand{\Lop}{\ensuremath{\mathcal L}}
\newcommand{\Mini}{\mathbb{M}}

\newcommand{\fa}{\boldsymbol{f}_\alpha}
\newcommand{\Sa}{{{}_{\,\alpha\!}\mathcal{S}}}
\newcommand{\Scortho}{{{}_{\,\pi/4\!}\mathcal{S}}}
\newcommand{\Sc}{{\mathcal{S}}}
\newcommand{\Score}{{\mathcal{S}_{core}}}
\newcommand{\Sext}{{\mathcal{S}_{ext}}}
\newcommand{\groupS}{{\mathscr{G}_\Sc}}  
\newcommand{\group}{{\mathscr{G}}}  
  
\newcommand{\groupsystem}{{\mathscr{G}_{system}}}  

\newcommand{\iequal}{{i=\pm1,\pm2}}
\newcommand{\ialpha}{{i,\alpha}}
\newcommand{\ipi}{{i,\pi/4}}
\newcommand{\jq}{{j,q}}
\newcommand{\ji}{{j,i}}
\newcommand{\althph}{{\alpha,\vartheta,\varphi}}
\newcommand{\althphp}{{\alpha,\vartheta,\varphi,\phi}}

\newcommand{\ialthphp}{{i,\alpha,\vartheta,\varphi,\phi}}
\newcommand{\onealpha}{{1,\alpha}}

\newcommand{\Stwo}{{\Sph^2}}
\newcommand{\Rfour}{{\Rbb^4}}
\newcommand{\Rthree}{{\Rbb^3}}
\newcommand{\Rtwo}{{\Rbb^2}}
\newcommand{\Rtwoplus}{{\Rbb^2_+}}

\newcommand{\vece}{\vec{e}}
\newcommand{\veceta}{\vec{\eta}}
\newcommand{\vecnu}{\vec{\nu}}
\newcommand{\vecmu}{\vec{\mu}}
\newcommand{\veceperp}{{\vec{e}^\perp}}
\newcommand{\vecep}{\vec{e}^{\,} {}'_{\!\ialthphp}}
\newcommand{\vecepji}{\vec{e}^{\,} {}'_{\!\ji}}

\newcommand{\vecepp}{\vec{e}^{\,} {}''_{\!i}}
\newcommand{\vecepm}{\vec{e}^{\,} {}''_{\!-i}}

\newcommand{\mmer}{{m_{mer}}}
\newcommand{\mpar}{{m_{par}}}
\newcommand{\mparp}{{m'_{par}}}
\newcommand{\Lmer}{{L_{mer}}}
\newcommand{\Lpar}{{L_{par}}}
\newcommand{\Ipar}{{I_{par}}}
\newcommand{\Hgluing}{{H_{gluing}}}
\newcommand{\Herror}{{H_{error}}}
\newcommand{\Hprescribed}{{H_{prescribed}}}

\begin{document}

\title[Doubling and Desingularization]{Doubling and Desingularization constructions for minimal surfaces}

\author[N.~Kapouleas]{Nikolaos~Kapouleas}
\address{Department of Mathematics, Brown University, Providence,
RI 02912} \email{nicos@math.brown.edu}





\begin{abstract}
In the first part of the paper we discuss the current status of
the application of the gluing methodology to
doubling and desingularization constructions
for minimal surfaces in Riemannian three-manifolds.
In particular a doubling construction for equatorial spheres in $\Sph^3(1)$ is announced.
Aspects of the current understanding of existence and uniqueness questions
for closed minimal embedded surfaces in $\Sph^3(1)$ are also discussed,
and some new uniqueness questions are proposed.
In the second part of the paper we discuss some of the ideas and provide an outline
for a general desingularization construction without imposed symmetries.

This paper is the author's contribution to the
volume in honor of Professor Richard M. Schoen's sixtieth birthday.
\end{abstract}

\maketitle

\section{Introduction}
\label{Sintrod}
\nopagebreak

\subsection*{Definitions and background}
$\phantom{ab}$
\nopagebreak

We start by defining doubling and desingularization.
Our definitions are consistent with the terminology introduced in
\cite{kapouleas:survey,kapouleas:finite,kapouleas:clifford}.

\addtocounter{equation}{1}
\begin{definition}
\label{Ddoubling}
We call a smooth minimal surface $M$,
a doubling of a minimal surface $\Sigma$ in a Riemannian manifold $\N$,
if we can write $M=\Sigma_1\cup\Sigma_2\cup \Sigma_{bridges}$
with $\Sigma_1$, $\Sigma_2$, and $\Sigma_{bridges}$ as follows:
\newline
(i).
$\Sigma_1$ and $\Sigma_2$ are graphs of two
functions (more generally sections of the normal bundle in $\N$)
$\phi_1$ and $\phi_2$ over $\Sigma\setminus\cup_{i=1}^N D_i$,
the given surface with $N$ small discs $D_i$ removed.
\newline
(ii).
$\Sigma_{bridges}=\cup_{i=1}^N B_i$
is the union of $N$ annuli $B_i$.
Each $\partial B_i$ is the union of the boundary circles
of $\Sigma_1$ and $\Sigma_2$ over $\partial D_i$.
\end{definition}

Note that there is no restriction on the dimensions of the minimal surfaces and $\N$.
In this paper however we concentrate on the case the surfaces are two-dimensional and $\N$
three-dimensional.
We expect the annuli $B_i$ in the above definition to resemble catenoidal bridges.
Given a doubling $M$ of $\Sigma$ as above, we can associate the set $L\subset\Sigma$ of centers of the discs
$D_i$.
Clearly $L$ is not precisely defined except in cases of high symmetry.
In constructing $M$ however we can consider $L$ as part of the data along with $\Sigma$.
$L$ then determines the number and approximate location of the discs $D_i$ and the bridges $B_i$.
Because of the fundamental role $L$ plays in our constructions we have the following.

\addtocounter{equation}{1}
\begin{definition}
\label{Dconfig}
We call the set $L\subset\Sigma$ the configuration of the doubling $M$.
\end{definition}

We mention now some
early constructions which can be considered as doublings:
Karcher, Pinkall, and Sterling \cite{KPS},
used Lawson's method \cite{L2} based on solving the Plateau problem for polygons,
to construct highly symmetric minimal surfaces
resembling doublings of the equatorial sphere in $\Sph^3(1)$.
Each of the surfaces has the symmetries of a Platonic solid
and the configuration $L$ consists of its vertices.
Because the Platonic solids are finitely many,
this provided only finitely many examples.
Pitts and Rubinstein have discussed \cite{PRu}
constructions by min-max methods for discrete families of minimal surfaces,
where the size of the catenoidal bridges used can be arbitrarily small,
and then the number of the bridges (and hence the genus)
tends to infinity,
while the surfaces tend to a limit varifold.
These constructions are also highly symmetric.
Some of the constructions resemble doubling constructions and have a limit varifold which is a minimal
surface counted with multiplicity two.
Wohlgemuth \cite{wohlg}
constructed minimal surfaces in $\Rthree$
which resemble two catenoids connected by a ``ring'' of catenoidal bridges.
These examples are also highly symmetric,
and in the limit as the size of the bridges tends to $0$, their number $N$
tends to infinity, and the minimal surfaces tend to a doubly covered catenoid.
The configuration $L$ is contained in the waist of the limit catenoid.

In the desingularization case the role of $\Sigma$ and $L$ is served
by a given minimal two-surface 
immersed in a three-dimensional Riemannian manifold $\N$
and a curve of intersection $\Cunder$ in $\N$.
The minimal surface is described
by a minimal immersion $\X:\W\to\N$
where $\W$ is an abstract surface which 
usually has many connected components.
It is unclear how to generalize successfully to higher dimensions.
As an example consider the desingularization of two coaxial catenoids 
intersecting along two circles.
$\W$ is then the abstract disjoint union of the two catenoids
and $\Cunder$ the union of two circles of intersection.
A more precise discussion of the given system of minimal surfaces is postponed until Section \ref{Sbuildingblocks}.
To avoid burdening the presentation too much we often assume that $\Cunder$ is embedded
although this is not needed for the theorem (see also \ref{CCunder}).
In this spirit
we often refer to $\Cunder$ as a curve without clarifying whether it is embedded or merely immersed.

\addtocounter{equation}{1}
\begin{definition}
\label{Ddesingularization}
Given $\W$ and $\Cunder$ as above,
we call a smooth minimal surface $M$ a desingularization of $\W$ along $\Cunder$,
if we can write $M=\underline{M}\cup\widehat{M}$ where 
$\underline{M}$ and $\widehat{M}$ are as follows:
\newline
(i).
$\widehat{M}=\Xhat'(\What)$,
where $\Xhat':\What\to\N$ is a perturbation of
the given minimal immersion $\Xhat:\What\to\N$,
where $\What$ is $\W$ modified by ``cutting'' along $\Cunder$ (see \ref{DWhat}).
Note that $\partial\What\setminus\partial\W$ covers $\Cunder$ four-to-one.
\newline
(ii).
$\underline{M}$ is contained in a tubular neighborhood of $\Cunder$
and is embedded when $\Cunder$ is.
Moreover
$\partial\underline{M}
=
\Xhat'
(\,
\partial\What\setminus\partial\W
\,)$.
\end{definition}

The term ``desingularization'' is motivated from the fact that the points of $\Cunder$
are singular in the sense that the tangent cone at them consists of more than one plane.
A neighborhood of $\Cunder$ is replaced by a smooth surface $\underline{M}$ which has no singular points.
It is expected that the price of removing the singular points is the many handles
included in $\underline{M}$
which makes the genus of $M$ much higher than the genus of $\W$.

We mention now some
early constructions which can be considered as desingularization constructions.
The famous surfaces $\xi_{m,k}$ with $k=1$ of Lawson \cite{L2} can be interpreted as desingularizations
of two orthogonal equatorial two-spheres in the three-sphere along a common equatorial circle.
Note that for $k>1$ the surfaces $\xi_{m,k}$ can be interpreted as desingularizations of $k+1$
spheres intersecting along a common circle, where in \ref{Ddesingularization}.i
we would have that
$\partial\What\setminus\partial\W$ is a $2(k+1)$-to-one covering of $\Cunder$ instead of four-to-one.
Such desingularizations where more than two surfaces intersect along a common curve is an exceptional occurrence.
A more generic case involves the construction of
certain minimal surfaces by variational or Enneper-Weierstrass methods
by Hoffman and Meeks \cite{HM2,HM3,HM5,HK}.
These surfaces are desingularizations of two intersecting coaxial catenoids or a catenoid and a plane.

\subsection*{Gluing constructions}
$\phantom{ab}$
\nopagebreak

In this paper we survey the use of gluing methods to carry out 
doubling and desingularization constructions.
To apply these methods given minimal surfaces are combined to provide
more complicated initial surfaces which are approximately minimal.
The initial surfaces are perturbed then to minimality by solving 
the appropriate partial differential equation.
These constructions share two important features with a number of other gluing constructions
\cite{kapouleas:annals,kapouleas:jdg,kapouleas:drops,kapouleas:wente,haskins:kapouleas:invent,haskins:kapouleas:survey}:
First, the manifolds or surfaces can be subdivided 
into regions which in some appropriate sense carry kernel,
and regions which resemble ``necks'' connecting the previous regions,
and which do not carry kernel.
Second, they are constructions in submanifold geometry, which leads to conditions
like the ``flexibility'' condition 
which do not appear in constructions for metrics.
R. Schoen in \cite{schoen} pioneered this kind of constructions in cases where the first feature holds
and also supervised the thesis of the author \cite{kapouleas:annals,kapouleas:jdg}.

The methodology was systematized and refined further in order to carry out
a challenging gluing construction for the Wente tori
\cite{kapouleas:wente:announce,kapouleas:wente}.
We have discussed this methodology in an earlier survey paper \cite{kapouleas:survey}
in which we reviewed the general method and how it is applied in various cases including
doubling and desingularization constructions.
This paper complements and extends \cite{kapouleas:survey}
in the case of doubling and desingularization constructions.
In particular we present progress in the case of doubling constructions,
discuss the proof in a general desingularizing construction announced in
\cite[Theorem F]{kapouleas:survey},
and discuss applications and related open questions and motivations,
including minimal surfaces in the round three-sphere (see Section \ref{Ssph}).
This paper has little overlap with \cite{kapouleas:survey} and is self-contained.

In a doubling construction we assume given
a minimal surface $\Sigma$ in a Riemannian manifold $\N$. 
Such a construction requires then the determination (implicitly or explicitly)
of a configuration $L$ (recall \ref{Dconfig}),
functions $\phi_1$ and $\phi_2$ as in \ref{Ddoubling},
and the size of each catenoidal bridge.
Doubling constructions by gluing methods are motivated by the expectation that as the number (locally)
of catenoidal bridges tends to infinity, their size will tend to zero, and the catenoidal bridges appropriately
blown-up tend to actual catenoids.

Reversing this, a gluing construction can be attempted where the initial surfaces are constructed 
in accordance with definition \ref{Ddoubling}, but fail to be minimal, and where each bridge
approximates an appropriately truncated small catenoid.
Balancing and matching considerations are expected to restrict the possible configurations $L$
and determine the functions and sizes in terms of $L$ (see Section \ref{Sdoubling}).
As usual in this kind of gluing construction,
one has to introduce appropriate parameters in the construction of the initial surfaces
and obtain in this way a whole family of initial surfaces.
One of the initial surfaces can then be perturbed to minimality.

Such a construction in full generality is not yet understood.
There has been satisfactory progress however in cases of high symmetry.
We will discuss the corresponding constructions in Section \ref{Sdoubling}.
We can distinguish two cases for these constructions.
In the first case $L$ is uniformly distributed on the whole surface.
In the second case $L$ is contained on a union of curves.
In the latter case there are some similarities with the desingularization constructions.

Desingularization constructions by gluing methods were
inspired by the observation by Hoffman and Meeks \cite{HM4}
that the singly periodic Scherk surfaces \cite{S,DHKW,Ni}
appear as blow-up limits of desingularizations at the curve of desingularization.
The simplest and most symmetric of the singly periodic Scherk surfaces are asymptotic to
two orthogonal planes and are given by the equation
\addtocounter{theorem}{1}
\begin{equation}
\label{EScherkortho}
\sinh x_1 \sinh x_2=\sin x_3.
\end{equation}
More generally the singly periodic Scherk surfaces form 
a one-parameter family of embedded minimal surfaces
$\Sa$,
where 
$\alpha\in(0,\pi/2)$.
$\Sa$
is asymptotic to four half-planes two of which form an angle $2\alpha$.
It is given \cite{Ni} by the equation 
\addtocounter{theorem}{1}
\begin{equation}
\label{EScherk}
\cos^2\!\alpha\,\cosh\frac {x_1}{\cos\!\alpha}
-
\sin^2\!\alpha\,\cosh\frac {x_2}{\sin\!\alpha}
=\cos x_3,
\end{equation}
and is clearly invariant under translation by $(0,0,2\pi)$.
In the literature the surfaces $\Sa$
are referred to as Scherk's fifth surfaces,
Scherk's singly periodic surfaces,
or Scherk-towers
\cite{Ni,HM4,DHKW,kapouleas:finite}, but
we will simply call them Scherk surfaces.

In a desingularization construction we assume given 
a minimal two-surface $\W$,
in a three-dimensional Riemannian manifold $\N$,
and a curve of intersection $\Cunder$,
as in \ref{Ddesingularization}.
(See Section \ref{Sbuildingblocks} for a more precise discussion.)
To carry out a gluing construction initial surfaces which are only approximately minimal,
but otherwise conform to the properties described in Definition \ref{Ddesingularization},
are constructed first.
The desingularizing part $\underline{M}$ of the initial surfaces should be carefully modeled
after the Scherk surfaces described above in \ref{EScherk}.
The number of handles introduced by each component of $\underline{M}$,
that is the number of handles used to desingularize each component of $\Cunder$,
is prescribed.
The gluing construction is expected to work when each of these numbers is large enough
depending on the given system of minimal surfaces.

We often impose a group of symmetries $\group\subset\groupsystem$ on the construction,
where $\groupsystem$ is defined as follows.
\addtocounter{equation}{1}
\begin{definition}
\label{Dgroupsystem}
For a given system of minimal surfaces as above,
we define the group of symmetries of the given system $\groupsystem$,
to be the group of isometries of the Riemannian manifold $\N$ preserving $\W$ and $\Cunder$.
\end{definition}
$\group$ has to be consistent with the symmetry group of the Scherk surfaces as well,
and therefore depends on the number of handles prescribed for the construction.

As in the case of doubling constructions,
general constructions with little symmetry imposed (small $\group$),
are much harder than highly symmetric constructions (appropriately large $\group$).
When $\group$ is small,
even the construction of the initial surfaces
involves solving Partial Differential Equations
and carefully estimating their solutions.
In Section \ref{Sdesingularization} we survey both highly symmetric and general constructions.
The second part of this paper presents
an outline of the construction and proof in a general case with small or trivial $\group$.

\addtocounter{equation}{1}
\begin{remark}
\label{Rcmc}
These constructions cannot be applied to constant mean curvature surfaces,
because as can be seen by inspection
the direction of the Gauss map of the surfaces constructed
cannot be chosen consistently with the direction of the Gauss map of the given surfaces.
\end{remark}

\subsection*{Motivation and applications}
$\phantom{ab}$
\nopagebreak

A motivation for general doubling or desingularization constructions is problem 88
in the list of open problems proposed by S.-T. Yau in 1982 \cite{yau:seminar}.
In this problem it is required to establish that there are infinitely many minimal
surfaces in any three-dimensional Riemannian manifold.
A general doubling construction would reduce this question to the existence of a single
minimal surface satisfying the appropriate necessary conditions.
Similarly our general desingularization Theorem \ref{Tmain} \cite{kapouleas:compact}
reduces this question to the existence of minimal surfaces appropriately intersecting
and satisfying the nondegeneracy conditions.
In either case then,
the doubling or desingularization construction would allow us to
conclude the existence of infinitely many minimal surfaces
by varying the number of catenoidal bridges or handles used in the construction.
Clearly to ensure the existence of one minimal surface satisfying the required conditions
in a large class of Riemannian manifolds,
it would help to have as general a doubling or desingularization theorem as possible.

Another potential application of the desingularization constructions is
to the Calabi-Yau problem for minimal surfaces in the embedded case.
Non-existence results have been proven in important cases by Colding and Minicozzi in \cite{CoMi}.
In other cases existence is expected,
and the recent successes of Martin, Meeks and their coauthors in the immersed case 
are very encouraging \cite{53,54}.
When existence is expected in the embedded case,
it seems that the main difficulty is to be able to remove the self-intersections by desingularizing.
Theorem \ref{Tmain} is not quite enough because it seems that in general one cannot
avoid triple intersection points \cite{pers}.
Allowing such points in the theorem in general,
requires finding and understanding minimal surfaces desingularizing three intersecting planes,
the way that Scherk surfaces desingularize two intersecting planes.
These surfaces could then be used as models in the desingularization constructions.

Another very interesting problem requires the classification and understanding of closed embedded
minimal surfaces in the round three-sphere, especially with an area bound (for example
the area of a few copies of the equatorial two-sphere) imposed.
At the moment there are some general theorems as for example \cite{choi:schoen,li:yau},
but few concrete examples known without using gluing constructions
\cite{L2,Ka3,PRu,KPS}.
As we discuss in the following sections, doubling \cite{kapouleas:clifford,kapouleas:equator}
and desingularization \cite{kapouleas:nested} constructions provide families of new interesting examples.
Actually it is not clear what nontrivial examples exist which cannot be interpreted as doublings or desingularizations.
Our constructions suggest interesting uniqueness questions also which we discuss in Section \ref{Ssph}.

Another source of many interesting doubling or desingularization constructions involves self-similar solutions for the mean curvature flow.
Self-shrinkers in particular can be considered as minimal surfaces in the Gaussian metric which is conformal to the Euclidean metric.
Certain highly symmetric desingularization constructions for self-shrinkers have been attempted by Nguyen \cite{nguyen1,nguyen3}.
Not all difficulties have been resolved at the moment \cite{nguyen5},
but there are many promising constructions and one expects many examples to be constructed in the future.
In the self-translating surfaces case, which can also be interpreted as minimal surfaces \cite{solitons},
Nguyen has carried out some highly symmetric desingularization constructions
\cite{nguyen2,nguyen4}.

Finally there are some interesting minimal surface
examples mentioned in \cite{kapouleas:survey},
which could be constructed if certain conditions for \ref{Tmain} are appropriately relaxed.
One such example requires desingularizing two intersecting catenoids
one of which is a translation of the other so that their waists are on the same plane but do not intersect.

\subsection*{Organization of the paper}
$\phantom{ab}$
\nopagebreak

This paper has two parts.
In the first part we review and announce the various results in the subject
and we also discuss various open questions.
In Section \ref{Sdoubling} we survey doubling constructions by gluing methods
and in Section \ref{Sdesingularization} desingularization constructions also by gluing methods.
In Section \ref{Ssph} we discuss closed embedded minimal surfaces in the round three-sphere.

In the second part of the paper we provide an outline of the construction and proof of Theorem \ref{Tmain}.
In Section \ref{Sbuildingblocks} we set the notation 
and discuss in some detail the given system of minimal surfaces and the family of singly periodic Scherk surfaces
which are the building blocks of the construction.
In Section \ref{Scentral} we present the construction of an initial surface.
In Section \ref{Sfamily} we discuss the modifications introduced to the construction of the remaining
initial surfaces needed in the construction and which form a smooth family of many parameters.
In the last Section \ref{Smain} we discuss some of the estimates and main ideas of the proof of \ref{Tmain}.

\subsection*{Notation and conventions}
$\phantom{ab}$
\nopagebreak

In this paper we adopt the convention that the mean curvature of the round two-sphere in Euclidean
three-space is $2$.
We have then that the mean curvature $H_f$ of the graph of a function $f$ over a surface of mean curvature $H$
is given by
\addtocounter{theorem}{1}
\begin{equation}
\label{EHf}
H_f=H+\Lop f + \mathcal{Q}_f,
\qquad
\text{where}
\qquad
\Lop:=\Delta+|A|^2+Ric(\nu,\nu),
\end{equation}
and $\mathcal{Q}_f$ is quadratic and higher order in $f$ and its derivatives
with coefficients involving the geometric invariants of the original surface,
$|A|^2$ is the square of the length of the second fundamental form $A$ of the original surface,
and $Ric(\nu,\nu)$ is the Ricci curvature of the ambient manifold evaluated on the unit normal of 
the original surface.

In this paper we use weighted H\"{o}lder norms.
The definition we use is given by
\addtocounter{theorem}{1}
\begin{equation}
\label{E:weightedHolder}
\|\phi: C^{k,\beta}(\Omega,g,f)\|:=
\sup_{x\in\Omega}\frac{\,\|\phi:C^{k,\beta}(\Omega\cap B_x,g)\|\,}{f(x)},
\end{equation}
where $\Omega$ is a domain inside a Riemannian manifold $(M,g)$,
$f$ is a weight function on $\Omega$,
$B_x$ is a geodesic ball centered at $x$ and of radius the minimum of
$1$ and half the injectivity radius at $x$.

We will be using extensively cut-off functions,
and for this reason we adopt the following notation:
We fix a smooth function $\psi:\Rbb\to[0,1]$ with the following properties:
\newline
(i).
$\psi$ is nondecreasing.
\newline
(ii).
$\psi\equiv1$ on $[1,\infty]$ and $\psi\equiv0$ on $(-\infty,-1]$.
\newline
(iii).
$\psi-\frac12$ is an odd function.
\newline
Given then $a,b\in \Rbb$ with $a\ne b$,
we define a smooth function $\psi[a,b]:\Rbb\to[0,1]$
by
\addtocounter{theorem}{1}
\begin{equation}
\label{Epsiab}
\psi[a,b]=\psi\circ L_{a,b},
\end{equation}
where $L_{a,b}:\Rbb\to\Rbb$ is the linear function defined by the requirements $L(a)=-3$ and $L(b)=3$.
Clearly then $\psi[a,b]$ satisfies the following:
\addtocounter{equation}{1}
\begin{lemma}
\label{Lpsiab}
(i).
$\psi[a,b]$ is weakly monotone.
\newline
(ii).
$\psi[a,b]=1$ on a neighborhood of $b$ and 
$\psi[a,b]=0$ on a neighborhood of $a$.
\newline
(iii).
$\psi[a,b]+\psi[b,a]=1$ on $\Rbb$.
\end{lemma}

\subsection*{Acknowledgments}
I would like to thank Rick Schoen for his constant interest and support.

\section{Doubling constructions}
\label{Sdoubling}
\nopagebreak

\subsection*{Balancing and heuristic arguments}
$\phantom{ab}$
\nopagebreak

Balancing considerations imply constraints to the existence of a doubling $M$ of a given
minimal surface $\Sigma$ (recall \ref{Ddoubling})
with a given configuration $L$ (recall \ref{Dconfig}).
(Balancing for minimal surfaces 
is based simply on the first variation formula \cite{KKS,Si}.
For a general discussion in the current context see \cite{kapouleas:survey}.)
In this subsection we describe heuristic arguments based on balancing considerations.
We do not intend at the moment to include these arguments in rigorous proofs because
they are based on various assumptions
whose validity we do not attempt to check.
Moreover in more general constructions than the ones we have currently,
we may need to study significant correction terms we currently ignore.
Our purpose here is to derive certain terms in balancing expressions,
demonstrate their importance,
obtain indispensable insights in setting up our gluing constructions,
and motivate the actual proofs which are based on precise calculations.
We do expect that in many cases the arguments that follow describe the full picture.

We start by describing a heuristic balancing argument \cite{kapouleas:clifford},
which suggests a necessary condition for the existence of doublings
and also that the size of the catenoidal bridges is determined by $L$.
We describe the argument in the case $L$ is large and fairly uniformly distributed on $\Sigma$.
Let $p\in L$.
Consider the catenoidal bridge centered at $p$.
In this discussion we refer to directions parallel or orthogonal to $\Sigma$ as horizontal or vertical
respectively.
Consider the region $\widetilde{\Omega}\subset M$ corresponding to half the catenoidal bridge
appropriately extended,
so that $\partial\widetilde{\Omega}=\partial_{waist}\cup\partial_1$,
where $\partial_{waist}$ is the waist of the bridge 
and $\partial_1$ is the graph of $\phi_1$ over
$\partial\Omega$,
where $\Omega\subset\Sigma$ is a neighborhood of $p$
such that the conormal to $\partial_1$ tangent to $M$ is approximately horizontal
and $\phi_1$ on $\partial_1$ approximately constant.
More precisely we require $\int_{\partial\Omega}\vec{\eta}(\phi_1)=0$,
where $\vec{\eta}$ is the outward unit normal to $\partial\Omega$ in $\Sigma$,
and that $\phi_1$ oscillates little compared to its average on $\partial\Omega$.
Such an $\Omega$ should exist because if we vary $\Omega$ from small
to larger the sign of the integral clearly starts positive and changes to negative.

Assuming that $\Omega$ and $\widetilde{\Omega}$ are small, 
we have that vertical translation is an approximate Killing field and we can use it to obtain
an approximate balancing formula for the vertical force through $\partial_{waist}$ and $\partial_1$.
Actually in the case when $M$ is the round sphere there are exact Killing fields which are perturbations of the 
vertical translation in the vicinity of $p$.
To calculate the force through $\partial_1$ we modify and smoothly extend $\phi_1$
so that it is defined on the whole $\Omega$,
while it is kept unchanged in a small neighborhood of $\partial\Omega$,
and oscillates little on $\Omega$.
We then use the balancing formula on the graph of $\phi_1$ to calculate the vertical force
through $\partial_1$ as an integral over $\Omega$.
Assuming that the nonlinear terms for the mean curvature of the graph are negligible (recall \ref{EHf}),
we conclude \cite[equation 1.1]{kapouleas:clifford}
\addtocounter{theorem}{1}
\begin{equation}
\label{Earea}
\text{Area}(\Omega)\,\, (|A|^2+Ric(\nu,\nu))\,\,\phi_1=2\pi\tau,
\end{equation}
where $\tau$ is the size (radius of the waist) of the catenoidal bridge.
Note that the values of the quantities involved are only approximately defined.
For example $|A|^2+Ric(\nu,\nu)$ in the formula should be
interpreted to mean some appropriate weighted average on $\Omega$.
When $\Omega$ is small $|A|^2+Ric(\nu,\nu)$ does not vary much on $\Omega$
and so its value in the formula is fairly well defined.
If $|A|^2+Ric(\nu,\nu)$ vanishes we should replace it with other higher order terms.

Since $\tau>0$ in \ref{Econdition},
the above heuristic argument suggests that a necessary condition for a doubling construction
is that the mean curvature of the parallel surfaces to $\Sigma$ points away from $\Sigma$,
or equivalently the surface area is unstable,
which amounts (unless $|A|^2+Ric(\nu,\nu)=0$) to
\addtocounter{theorem}{1}
\begin{equation}
\label{Econdition}
|A|^2+Ric(\nu,\nu) >0 \qquad\text{ on }\quad L \subset \Sigma.
\end{equation}
This condition ensures that the vertical components of the force through 
the components of 
$\partial\widetilde{\Omega}=\partial_{waist}\cup\partial_1$
point in opposite directions and so can cancel each other.

Another way to heuristically justify \ref{Econdition} as a requirement for a gluing construction
is as follows.
When an initial surface is constructed mean curvature is created by the extra bending on the catenoidal
bridge to attach it to a parallel surface to $\Sigma$, and this parallel surface
has mean curvature as well.
In order to be able to correct the initial surface we expect that these two contributions to the mean curvature
combine to make the mean curvature $L^2$-orthogonal to the constants.
\ref{Econdition} ensures that the signs of these two contributions are opposite,
making orthogonality possible.

By assuming now that the catenoidal bridge follows fairly closely
a trunca\-ted ca\-te\-no\-id,
we conclude that 
\addtocounter{theorem}{1}
\begin{equation}
\label{Eheight}
\phi_1\sim\tau\log\frac{d}{\tau}
\qquad
\text{ on }
\partial_1,
\end{equation}
where $d$ indicates the distance of the points of $\partial\Omega$ from $p$
and is therefore only defined up to some factor $C$.
Using this in \ref{Earea} we conclude that 
\addtocounter{theorem}{1}
\begin{equation}
\label{Etau}
\tau\sim \,\, d  \,\,    e^{-\frac{2\pi}{ \text{Area}(\Omega)\,\, (|A|^2+Ric(\nu,\nu))  }},
\end{equation}
where both $d$ and $\tau$ are determined only up to a certain constant.
This indicates that $L$ determines the size of the catenoidal bridges.

In certain doublings,
$L$ is not uniformly distributed on $\Sigma$,
and then it is not possible to find a small (in diameter) $\Omega$ so that
$\int_{\partial\Omega}\vec{\eta}(\phi_1)=0$.
If we insist that $\Omega$ is small but allow nonvanishing $\int_{\partial\Omega}\vec{\eta}(\phi_1)$,
we can establish by the same argument as before that
\ref{Earea} is modified to 
\addtocounter{theorem}{1}
\begin{equation}
\label{Evertical}
\left(
\text{Area}(\Omega)\,\, (|A|^2+Ric(\nu,\nu))\,\,\phi_1
+\int_{\partial\Omega}\vec{\eta}(\phi_1)
\right)
\sim 2\pi\tau,
\end{equation}
where we ignore nonlinear terms in the derivatives of $\phi$.

Finally we discuss the horizontal components of the force through $\partial_1$.
As before we assume that a neighborhood of $\partial_1$ in $M$ is the graph over a neighborhood of
$\partial\Omega$ by a function $\phi_1$.
Moreover to simplify the argument we assume that $\phi_1$ is constant on $\partial\Omega$,
and that there exists a Killing field $\vecK$ in the vicinity of $p$,
which is exactly horizontal and approximates a horizontal translation.
The force with respect to $\vecK$ is defined to be
$$
F_\vecK=\int_{\partial_1}\vec{\eta}_1\cdot\vecK,
$$
where $\vec{\eta}_1$ is the outward conormal to $\partial_1$ tangent to $M$.
Moreover by the balancing formula and our assumptions we have
$$
\int_{\partial_1}\vec{\eta}_{hor}\cdot\vecK=0,
$$
where $\vec{\eta}_{hor}$ is the horizontal outward unit normal to $\partial_1$.
By subtracting and ignoring cubic and higher order terms in $\phi_1$ and its derivatives,
we conclude that
\addtocounter{theorem}{1}
\begin{equation}
\label{Ehorizontal}
F_\vecK=-\frac12\int_{\partial_1}
(\vec{\eta}(\phi_1))^2
\,\,\vec{\eta}\cdot\vecK,
\end{equation}
where $\vec\eta$ is the outward unit normal to $\partial\Omega$ on $\Sigma$.

\subsection*{Doubling the Clifford torus with a square lattice configuration \cite{kapouleas:clifford}}
$\phantom{ab}$
\nopagebreak

We outline now a doubling construction for the Clifford torus with maximal symmetry.
Let $\Sigma=\Sph^2(\frac1{\sqrt2})\times\Sph^2(\frac1{\sqrt2})\subset\Sph^3(1)$
be the Clifford torus and $L\subset\Sigma$ an $m\times m$ square lattice on the Clifford torus,
where we assume $m$ to be large.
The symmetry group $\group$ imposed on the construction consists of the isometries of $\Sph^3(1)$
which fix $L$ as a set.
$\group$ is quite large and this simplifies the construction because of the following:
\addtocounter{equation}{1}
\begin{lemma}
\label{Lgroup}
(i). 
$\group(L)=L$ by the definition of $\group$
and $\group(\Sigma)=\Sigma$.
\newline
(ii).
If $\Sigma'$ is a small perturbation 
of $\Sigma$ with 
$\group(\Sigma')=\Sigma'$, then $\Sigma'=\Sigma$.
\newline
(iii). $\group$ acts transitively on $L$.
\newline
(iv). 
If $L'$ is a small perturbation in $\Sph^3(1)$ of $L$ with 
$\group(L')=L'$, then $L'=L$.
\newline
(v). There is a symmetry in the stabilizer $\group_p\subset \group$ of $p$ in $\group$
which exchanges the two sides of $\Sigma$,
but there are no such symmetries which fix $\Sigma$ or $L$ pointwise.
\end{lemma}

This lemma follows easily by using a coordinate system on the unit sphere in which the metric takes the form
\addtocounter{theorem}{1}
\begin{equation}
\label{EPhig}
    g = (1+\sin2  z)\,d  x^2 + (1-\sin2  z)\,d  y^2 + d  z^2,
\end{equation}
and in which the Clifford torus corresponds to $\{z=0\}$ (see \cite{kapouleas:clifford} for details).

We construct now the initial surfaces which depend smoothly on a parameter $\tau$.
Given a small $\tau>0$ we define a smooth, embedded, connected,
and closed, surface $M_\tau$
by
\addtocounter{theorem}{1}
\begin{equation}
\label{EMtau}
M_\tau=\Sigma_1\cup\Sigma_2\cup \Sigma_{bridges},
\end{equation}
where we have the following:
\newline
(a). $\Sigma_1$ and $\Sigma_2$ are parallel copies 
of $\Sigma\setminus\cup_{p\in L}D_p$
at heights $\pm\tau\log\frac d \tau$,
where $D_p$ is a geodesic disc in $\Sigma$ of radius $2d$ and center $p$,
and $d$ is a constant smaller than a fourth of the step of the lattice and of the same order,
for example we can take $d=1/m$.
\newline
(b). $\Sigma_{bridges}=\cup_{p\in L} B_p$,
where each $B_p$ is a truncated catenoid centered at $p$,
with a region of transition close to its boundary
to ensure smooth attachment to $\Sigma_1$ and $\Sigma_2$.
To describe the catenoids one can use
the coordinates in \ref{EPhig}.
This way although the bridges in $\Sigma_{bridges}$
are not minimal, they have small mean curvature in the appropriate sense.

$\tau$ is allowed to vary around the value determined by \ref{Etau}.
More precisely we allow values which are multiples up to a uniform (independent of $m$) factor
of the value
\addtocounter{theorem}{1}
\begin{equation}
\label{Etb}
  \tb := m^{-1}e^{-m^2/4\pi}.
\end{equation}
We have then the following.
\addtocounter{equation}{1}
\begin{theorem}[\cite{kapouleas:clifford}]
\label{Tclifford}
If $m$ is large enough, then there is a $\tau$ in the above range,
such that $M_\tau$ can be perturbed to a minimal surface $\Sigma_m$
which is smooth, embedded, connected and closed.
$\Sigma_m$ is a doubling of the Clifford torus $\Sigma$ in the sense of \ref{Ddoubling},
with configuration $L$ in the sense of \ref{Dconfig},
where $L\subset \Sigma$ is the $m\times m$ square lattice discussed above.
Moreover as $m\to\infty$, $\Sigma_m$ tends as a varifold to $\Sigma$ with multiplicity two. 
\end{theorem}

The general idea for the proof of this theorem is as follows and as mentioned already is in accordance
with the methodology developed in \cite{kapouleas:wente:announce,kapouleas:wente,kapouleas:imc},
discussed in \cite{kapouleas:survey},
and used in various other constructions
\cite{kapouleas:finite,haskins:kapouleas:invent,haskins:kapouleas:survey}:
On each $M_\tau$ we solve a Partial Differential Equation to find a function
whose graph has mean curvature in the (extended) substitute kernel $\skernel$,
which in general is a function linear space corresponding to eigenfunctions of small eigenvalue
and certain low harmonics on certain meridians.
As we discuss below, it turns out that in our case $\skernel$ is one-dimensional.
Balancing (using a precise calculation) implies the existence of a $\tau$
for which an integral of the mean curvature of the graph vanishes.
By the one-dimensionality of $\skernel$ this is enough to ensure the minimality of the corresponding graph.

To understand $\skernel$ we use 
\addtocounter{theorem}{1}
\begin{equation}
\label{Eh}
h:=\frac{|A|^2+m^2}{2}g,
\end{equation}
a metric conformal to the induced metric $g$ on $M_\tau$.
The relevant linear operator then is given by 
\addtocounter{theorem}{1}
\begin{equation}
\label{ELh}
\Lh:=\Delta_h+2\frac{|A|^2+2}{|A|^2+m^2}.
\end{equation}

We concentrate our attention now to $M_{\tau,p}$,
the set of points of $M$ which are closer to $p\in L$ than any other point in $L$.
Clearly $\group M_{\tau,p}=M$.
It turns out that $M_{\tau,p}$ equipped with $h$ tends as $m\to\infty$
to the union of two flat squares and a round unit sphere.
The linear operator $\Lh$ tends to $\Delta+2$ on the unit sphere and 
the flat Laplacian with Neumann boundary conditions on the squares.
The kernel then in the limit is five-dimensional with corresponding eigenfunctions
the constants on the squares and the first harmonics on the sphere.
The symmetries of the construction however kill the first harmonics on the sphere
and identify the two squares, reducing this way the dimension of the kernel to one,
with corresponding eigenfunction a constant on the squares.
The approximate kernel on $M_\tau$ is then one-dimensional and it turns out that
we do not need any extra extended substitute kernel to ensure appropriate decay.
This implies that the extended substitute kernel $\skernel$ is one-dimensional.

\subsection*{Doubling the Clifford torus with a rectangular lattice configuration \cite{wiygul}}
$\phantom{ab}$
\nopagebreak

The construction described above is the ``most'' symmetric construction one can hope for
as demonstrated by Lemma \ref{Lgroup}.
It is natural to try to understand less symmetric constructions 
by gradually relaxing the symmetry assumptions.
The simplest first generalization is to replace the square lattice $L$
on the Clifford torus with
rectangular lattices $k_1 m\times k_2 m$ 
where $k_1,k_2$ are relatively prime and 
$m$ is large enough in terms of $k_1,k_2$.
In such a configuration conditions (ii), (iv), and (v) in \ref{Lgroup} fail,
as is evident from \ref{EPhig}.

The other conditions are still satisfied however and condition (iv) holds if it is weakened by
the assumption $L'\subset\Sigma$.
This means that there is enough symmetry so that the catenoidal bridges can be identified with each other.
Eventually one has to deal with a three-dimensional kernel:
There is no symmetry to identify the two copies of the Clifford torus,
and so the two squares in the fundamental domain carry two different (constant)
eigenfunctions in the limit as $\tau\to0$.
First harmonics corresponding to vertical translations on the catenoidal region
(which tends to a round sphere in the $h$ metric) are also allowed,
but not the first harmonics corresponding to horizontal translations.
The three-dimensional kernel means that two more parameters have to be introduced
in the construction of the family of initial surfaces:
One parameter corresponds to moving the whole initial surface vertically
so that the heights of two copies of the Clifford torus on its two sides
are not the same anymore.
The other parameter controls a change of size between the two halves of the catenoidal bridge.

\subsection*{Doubling the equatorial two-sphere in $S^3(1)$ \cite{kapouleas:equator}}
$\phantom{ab}$
\nopagebreak

We consider now the case $\Sigma=\Sph^2(1)\subset\Sph^3(1)$,
an equatorial two-sphere inside a unit three-sphere.
The main difficulty in doubling $\Sigma$ is that there are only finitely many symmetry groups
(related to the symmetry groups of the Platonic solids)
for which \ref{Lgroup}.iv applies.
Since our gluing approach works only when the number of points in $L$ is large enough,
we cannot avoid working with symmetry groups which allow the points of $L$ to ``slide''.
This means that the corresponding catenoidal bridges will have to find their position
by ``horizontal balancing'', and therefore horizontal forces and interactions have
to be understood.

In order to simplify the constructions as much as possible we should still impose the maximal possible symmetry.
The configurations we consider are as follows and were proposed some time ago by Hermann Karcher.
We assume that an equator circle and poles have been chosen on $\Sigma$,
and therefore the meridian semicircles and parallel circles of constant latitude have been determined.
We assume we are given $\mmer,\mpar\in\Nbb$.
$L$ then consists of the $\mmer\mpar$ intersection points
of $\mmer$ meridians with $\mpar$ parallels.
Because we require maximal symmetry the $\mmer$ meridians are arranged at equal angles around $\Sigma$
and successive meridians make an angle $2\pi/\mmer$.
The $\mpar$ parallels include the equator circle if $\mpar$ is odd.
The remaining parallels are positioned symmetrically on the two sides of the equator circle.
We number the parallels by increasing distance from the equator.
To facilitate reference we define $\Ipar:=\{\pm1,\pm2,...,\pm\mparp\}$ if $\mpar$ is even,
and $\Ipar:=\{0,\pm1,\pm2,...,\pm\mparp\}$ if $\mpar$ is odd,
where $\mparp$ is the integer part of $\mpar/2$.

The latitudes of the parallels are free to vary,
so to determine the configuration $L$ we specify
either the latitudes $x_i$ or their lengths $\ell_i$,
where $i\in\Ipar$.
We have then
\addtocounter{theorem}{1}
\begin{equation}
\label{Elengths}
\begin{gathered}
x_{-i}=-x_i,
\qquad
\ell_{-i}=\ell_i
=2\pi\cos x_i,
\\
0<x_1<x_2< ... < x_\mparp <\pi/2,
\qquad
2\pi>\ell_1>\ell_2> ... > \ell_\mparp >0.
\end{gathered}
\end{equation}
To help with the notation we use $\Lmer$ and $\Lpar=\Lpar(\{x_i\})$ to denote
the union of the meridians and parallels in consideration respectively.
We have then $L=\Lmer\cap\Lpar$.

The symmetry group $\group$ imposed on the construction consists as in the Clifford torus case
of the isometries of $\Sph^3(1)$ fixing $L$ as a set.
$\group$ does not depend on the latitudes of the parallels and the lengths $\ell_i$,
actually it depends only on $\mmer$.
We enumerate its properties in the following.

\addtocounter{equation}{1}
\begin{lemma}
\label{Lgroup2}
(i). 
$\group(L)=L$ by the definition of $\group$
and $\group(\Sigma)=\Sigma$.
\newline
(ii).
If $\Sigma'$ is a small perturbation 
of $\Sigma$ with 
$\group(\Sigma')=\Sigma'$, then $\Sigma'=\Sigma$.
\newline
(iii). $\group$ acts transitively on the meridians
and there is a symmetry which exchanges the parallels of the same length.
\newline
(iv). 
If $L'$ is a small perturbation in $\Sph^3(1)$ of $L$ with 
$\group(L')=L'$, then $L'=\Lmer\cap\Lpar(\{x'_i\})$,
where the $x'_i$'s are small perturbations of the $x_i$'s for which
$L=\Lmer\cap\Lpar(\{x_i\})$.
\newline
(v). 
There is a symmetry which exchanges the two sides of $\Sigma$ and keeps
$\Sigma$ fixed pointwise.
\end{lemma}

Note that \ref{Lgroup2}.v is stronger than \ref{Lgroup}.v.
This relates to the fact that the equator two-sphere is totally geodesic but the Clifford torus is not.
On the other hand,
as we have already mentioned,
\ref{Lgroup2}.iv is much weaker than \ref{Lgroup}.iv and allows ``sliding'' of the necks.
\ref{Lgroup2}.iii is also weaker than \ref{Lgroup}.iii.
Because of these differences the current construction is harder than the construction in \cite{kapouleas:clifford}.

As before the initial surfaces consist of catenoidal bridges and graphs of functions on $\Sigma$ minus small discs.
The functions we use are not constant as in the previous construction.
The determination of the latitudes of the parallel circles is part of the determination of the functions.
To describe the process we have the following.

\addtocounter{equation}{1}
\begin{lemma}
\label{Lbalanced}
Given $\mpar$ and $\mmer$ as above with $\mmer/\mpar$ large and $\mpar\ge2$,
there is 
a unique continuous function $\phi>0$ on $\Sigma$ 
and a unique choice of $x_i$'s as in \ref{Elengths},
such that the following hold.
\newline
(i). $\phi$ is rotationally invariant and symmetric under reflection with respect to the equator circle.
Equivalently it depends only on $|x|$ where $x$ is the latitude.
\newline
(ii). It is smooth on $\Sigma\setminus\Lpar$ where it satisfies the linearized equation for a minimal graph
$$
(\Delta+2)\phi=0.
$$
Note that because of the rotational invariance this amounts to a second order ODE.
Note also that $\phi$ is then piecewise smooth with jump discontinuities on its derivatives on $\Lpar$.
\newline
(iii).
We define $\tau_i>0$ ($i\in\Ipar$) by 
$\phi(x_i)=\tau_i\log\frac{\ell_i}{\mmer\tau_i}$.
We require then (vertical balancing)
$$
2\pi\tau_i=
\frac{\ell_i}{\mmer}
\left(
\left. \frac{\partial\phi}{\partial x}\right|_{x=x_i+}
-
\left. \frac{\partial\phi}{\partial x}\right|_{x=x_i-}
\right),
$$
where ${x=x_i\pm}$ denotes the one-sided derivatives at $x=x_i$ from right and left respectively.
\newline
(iv).
On $\Lpar$ we require (horizontal balancing)
$$
\left. \frac{\partial\phi}{\partial x}\right|_{x=x_i+}
+
\left. \frac{\partial\phi}{\partial x}\right|_{x=x_i-}
=0,
$$
\end{lemma}

The $\tau_i$'s determine the size of the catenoidal bridges we use in the construction of the initial surfaces.
\ref{Lbalanced}.iii is motivated then from \ref{Evertical}.
\ref{Lbalanced}.iv is motivated from an appropriate modification of
\ref{Ehorizontal}.
To keep the presentation simple
we outline the (elementary) proof of \ref{Lbalanced}
only in the simplest case where $\mpar=2$.
We assume the existence of $\phi$ as required.
In this case $\Sigma\setminus\Lpar=\Omega_0\cup\Omega_1\cup\Omega_{-1}$,
where $\Omega_0$, $\Omega_1$, and $\Omega_{-1}$ are connected domains,
$\Omega_0$ contains the equator circle,
and $\Omega_1$ is contains the North pole
(where the latitude $x=\pi/2$).

By uniqueness of solutions of the ODE we have
\addtocounter{theorem}{1}
\begin{equation}
\label{EOO}
\phi=\phi(\pi/2)\sin x
\text{ on }
\Omega_1
\qquad
\text{ and }
\qquad
\phi=\phi(0)\varphi
\text{ on }
\Omega_0,
\end{equation}
where $\varphi$ is defined on $\Sigma$ by requiring the following.
$\varphi$ satisfies the conditions in \ref{Lbalanced}.i,
and also the ODE implied by $(\Delta+2)\varphi=0$,
with initial conditions $\varphi(0)=1$
and $\frac{\partial\varphi}{\partial x}(0)=0$.
Observe that the smallest eigenvalue of $\Delta+2$
on $U(t):=\{|x|\le t\}\subset\Sigma$ with Dirichlet boundary data,
varies from very large for small $t$,
to $-2$ for $t=\pi/2$.
We conclude there is a unique $t=t_0$,
for which the smallest eigenvalue vanishes.
By the uniqueness for the ODE we conclude that $\varphi(t_0)=0$,
and $\varphi(x)>0$ for $x\in(-t_0,t_0)$.
By integrating then $(\Delta+2)\varphi=0$ over $U(x)$
and using the divergence theorem we conclude
that for $x\in(0,t_0)$
\addtocounter{theorem}{1}
\begin{equation}
\label{EOder}
\frac{\partial\varphi}{\partial x}(x)=-\frac1{\cos(x)}\int_{U(x)}\varphi   < 0.
\end{equation}

Since we are assuming that $\phi$ satisfies the conditions of the lemma
we have $\phi>0$ which by the above implies that $x_1<t_0$.
Moreover by using \ref{EOO} and \ref{EOder} and canceling $\phi(\pi/2)$ and $\phi(0)$
we have 
$$
\frac1{\phi(x_1)}
\left(
\left. \frac{\partial\phi}{\partial x}\right|_{x=x_1+}
+
\left. \frac{\partial\phi}{\partial x}\right|_{x=x_1-}
\right)
=
\cot(x_1)
-
\frac1{\varphi(x_1)\, \cos(x_1)}\,\int_{U(x_1)}\varphi.
$$
The right-hand side as a function of $x_1\in(0,t_0)$ has the following properties:
It is strictly decreasing because $\cot$, $\cos$, and $\varphi$ are strictly decreasing and $\int_{U(x)}\varphi$ is strictly increasing.
It tends to $\infty$ as $x\to0+$,
and tends to $-\infty$ as $x\to t_0-$.
It has therefore a unique root which is $x_1$ by  \ref{Lbalanced}.iv.
We have therefore determined $x_1$ uniquely and $\phi$ up to a constant factor.
\ref{Lbalanced}.iii determines then the factor uniquely.
This completes the uniqueness part of the proof.
Existence follows by checking that the conditions are satisfied once we define $\phi$
by using the expressions above.

\addtocounter{equation}{1}
\begin{remark}
\label{Rmerone}
Note that the lemma is not true in the case $\mpar=1$.
In this case $\Lpar$ is simply the equator circle.
By uniqueness of the ODE we have then
$\phi=\phi(\pi/2)\sin x$ on the Northern hemisphere,
which implies that $\phi=0$ on $\Lpar$,
which contradicts \ref{Lbalanced}.iii.
\end{remark}

To continue with the construction we have to ``unbalance'' $\phi$,
as usually required by the methodology we employ.
To this effect we can extend \ref{Lbalanced} to apply to the case where
\ref{Lbalanced}.iii and \ref{Lbalanced}.iv
have been modified by addition in the equations of terms controlling ``unbalancing''.
We obtain families of $\phi$'s which we use then to construct families of initial surfaces as follows.

To construct the initial surfaces we first ``smooth out'' $\phi$ near $\Lpar$,
to obtain a smooth $\widetilde{\phi}$ which agrees with $\phi$ exactly
except on a small neighborhood of $\Lpar$.
We then define similarly to \ref{EMtau} 
a smooth, embedded, connected, and closed, surface $M_\zeta$,
where $\zeta$ are the values of the parameters controlling unbalancing
(when $\zeta=0$ $\phi$ is given by \ref{Lbalanced}),
by
\addtocounter{theorem}{1}
\begin{equation}
\label{EMzeta}
M_\zeta=\Sigma_1\cup\Sigma_2\cup \Sigma_{bridges},
\end{equation}
where we have the following:
\newline
(a). $\Sigma_1$ and $\Sigma_2$ are the graphs of $\pm\widetilde{\phi}$
over $\Sigma\setminus\cup_{p\in L}D_p$,
where $D_p$ is a geodesic disc in $\Sigma$ of radius $2d$ and center $p$,
and we take $d=\frac{\ell_i}{4\mmer}$ when $p$ is on a parallel of length $\ell_i$.
\newline
(b). $\Sigma_{bridges}=\cup_{p\in L} B_p$,
where each $B_p$ is a truncated catenoid centered at $p$,
with a region of transition close to its boundary
to ensure smooth attachment to $\Sigma_1$ and $\Sigma_2$.
To describe the catenoids we use Fermi coordinates around $\Sigma$ and geodesic
coordinates around $p$ in $\Sigma$.

\addtocounter{equation}{1}
\begin{remark}
\label{Rdoublingbyscherk}
An alternative construction of the initial surfaces would be to place
singly periodic Scherk surfaces of the appropriate scale and angle $\alpha$ (recall \ref{EScherk})
and fuse their wings with the graphs of $\pm\phi$.
This would make this case of doubling resemble the desingularization constructions although
the main difference is that here the angle $\alpha\to0$ as the genus tends to $\infty$.
There do not seem to be any major advantages in this approach,
and we have chosen the catenoid approach which provides uniformity of presentation with the other
doubling constructions.
\end{remark}

\addtocounter{equation}{1}
\begin{theorem}[\mbox{\cite{kapouleas:equator}}]
\label{Tequator}
If $\mpar$ and $\mmer$ are as above, with $\mpar\ge2$ and $\mmer$ large enough in terms of $\mpar$,
then one of the initial surfaces $M_\zeta$ described above
can be perturbed to a minimal surface $\Sigma_{\mpar,\mmer}$
which is smooth, embedded, connected and closed.
Moreover $\Sigma_{\mpar,\mmer}$ is a doubling of the equatorial two-sphere $\Sigma$
in the sense of \ref{Ddoubling}
and tends as a varifold to $\Sigma$ covered twice as ${\mmer}\to\infty$.
\end{theorem}

\addtocounter{equation}{1}
\begin{remark}
\label{Requatorcases}
With further work we hope to strengthen \ref{Tequator}
to apply for any $\mpar\ge2$, $\mmer\ge3$, with the product
$\mpar\mmer$ large enough in terms of an absolute constant.
This would include then the case where $L$ is (roughly) uniformly distributed on $\Sigma$
and the case where $L$ concentrates on (at least three) meridians.
\end{remark}

\subsection*{More doubling constructions and open questions}
$\phantom{ab}$
\nopagebreak

The doubling of the Clifford torus with $L$ a rectangular lattice where conditions (ii) and (v) in \ref{Lgroup} fail,
and (iv) is weakened,
and the doubling of the equatorial sphere we discussed above where conditions (iii) and (iv) are substantially weakened,
can serve as prototypes of highly symmetric doubling constructions where the symmetry is weaker than in the case of \cite{kapouleas:clifford}
and allows the difficulties we discussed to emerge.
Successfully dealing with these difficulties
seems to provide a framework for dealing with a large class of highly symmetric doubling constructions.
For example such doubling constructions seem to hold great promise in constructing self-similar solutions of the mean curvature flow.
Another test in gradually expanding the applicability of this approach is to understand
situations where conditions (iii) and (iv) in \ref{Lgroup} are weakened further.
Such configurations $L$ would arise when the symmetry group is kept the same in the constructions
for the Clifford torus or the equatorial sphere discussed above,
but more catenoidal bridges are introduced
in a way that their number per fundamental domain remains finite.

As for desingularization constructions,
general constructions where no symmetry can be imposed seem harder.
Unlike in the desingularization case no such doubling constructions are understood currently.
The main difficulty seems to involve determining the configuration $L$ (that is the positions of the bridges)
so they can balance appropriately under the horizontal forces.

\section{Desingularization constructions}
\label{Sdesingularization}
\nopagebreak

\subsection*{Highly symmetric constructions}
$\phantom{ab}$
\nopagebreak

Recall that in a desingularization construction we assume given 
a minimal two-surface $\W$,
in a three-dimensional Riemannian manifold $\N$,
and a curve of intersection $\Cunder$,
as in \ref{Ddesingularization}.
It is often the case that $\groupsystem$ (defined in \ref{Dgroupsystem})
is large and acts transitively on each component of $\Cunder$.
One can then impose a group of symmetries $\group\subset\groupsystem$ which depends on the number
of handles prescribed and has the following feature:
The number of handles introduced by the construction in a fundamental region of $\group$ is finite and fixed
independently of their total number
or (equivalently) their size.
We refer to such constructions as ``highly symmetric'' constructions.
As we will see later such constructions avoid many difficulties present in general constructions with little or no symmetry,
where the number of handles per fundamental region of $\group$ tends to $\infty$ as the size of the handles tends to $0$.

An early highly symmetric construction 
was carried out by M. Traizet \cite{T2}.
In this construction many of the usual difficulties present in a highly symmetric construction are avoided
since the system being desingularized is a union of intersecting planes.
The Scherk surfaces used decay to their asymptotic planes.

A typical highly symmetric construction was carried out independently by the author \cite{kapouleas:finite}.
In this construction the system of minimal surfaces being desingularized is a finite collection
of coaxial catenoids and planes intersecting along round circles.
The initial surfaces have to be constructed carefully so that the mean curvature decays exponentially away
from the circles of intersection.
Another important feature of this construction is that the surfaces being desingularized are complete but not compact.
Apart from this, this construction has many similarities with the special case of \ref{Tmain} where $\groupsystem$
acts transitively on each component of $\Cunder$.
In this paper we do not discuss the construction of \cite{kapouleas:finite} further.
We refer the reader instead to the original paper \cite{kapouleas:finite},
and for a  more general discussion to \cite{kapouleas:survey}.

Finally we remark that as mentioned in the introduction
interesting highly symmetric desingularization constructions have been attempted
for self-similar solutions of the mean curvature flow \cite{nguyen1,nguyen2,nguyen3,nguyen4}.

\subsection*{A general theorem for constructions with little or no symmetry}
$\phantom{ab}$
\nopagebreak

In the second part of this paper we outline the proof of the
following theorem which we announced in detail in \cite[Theorem F]{kapouleas:survey}.
For a general discussion and comparison with other gluing constructions
we also refer to \cite{kapouleas:survey}.

\addtocounter{equation}{1}
\begin{theorem}[\cite{kapouleas:compact}]
\label{Tmain}
We assume given
a minimal two-surface $\W$,
in a three-dimen\-si\-o\-nal Riemannian manifold $\N$,
and a curve of intersection $\Cunder$.
We make the following convenient assumptions:
\newline
(a). $\W$ is compact, perhaps with boundary.
$\Cunder$ is then compact also.
\newline
(b).
There are no points of triple intersection.
\newline
(c).
$\Cunder$ is a curve of transverse intersection.
(This with (b) implies that $\Cunder$ is smooth.)
\newline
(d). $\Cunder$ does not intersect the boundary of $\W$ and therefore by (b) and (c)
it is a disjoint union of smooth circles.

If we prescribe then large enough numbers of handles in desingularizing each component of $\Cunder$,
families of initial surfaces can be constructed so that one of the initial surfaces can be perturbed to be a desingularization
of $\W$ along $\Cunder$ in the sense of \ref{Ddesingularization},
provided that the following non-degeneracy assumptions hold:
\newline
(i).
The kernel for the linearized operator
$\Lop =\Delta+|A|^2+Ric(\nu,\nu)$ on $\W$,
with Dirichlet conditions on $\partial\W$,
is trivial (unbalancing condition).
\newline
(ii).
The kernel for the linearized operator
$\Lop$ on $\What$ (recall \ref{Ddesingularization}.i),
with Dirichlet conditions on $\partial\What$,
is trivial (flexibility condition).

The desingularizations moreover tend as varifolds
to $\W$ as the number of handles desingularizing each component tends
to $\infty$.
\end{theorem}

Note that conditions (a-d) above provide the simplest setting for a general desingularization theorem.
In the simplest version of the theorem which we discuss later we take the phrase
``large enough numbers of handles in desingularizing each component of $\Cunder$''
to mean that for some given large constant $\cunder_0$ there is a large $m_0$,
depending only on the given system of minimal surfaces and $\cunder_0$,
such that the construction works if for each component $\Cunder_j$ of $\Cunder$,
the number of handles used is $m_j/2$, where 
each $m_j>m_0$ and each ratio $m_j/m_{j'}$ is bounded by $\cunder_0$ (see \ref{Cmj} also).
Note also that the $m_j$ have to be either integers or half-integers depending on the topology.
In the simplest case when there are two one-sided distinct components of $\W$,
intersecting through the component of $\Cunder$ in consideration,
which is one-sided inside them,
clearly $m_j$ has to be integer.

Condition (ii) is called the ``flexibility condition'' because it ensures that a small
perturbation of $\Cunder$
provides boundary data for perturbing $\What$ so that it remains minimal,
while $\Cunder$, instead of (locally) belonging to two smooth surfaces intersecting through it,
becomes (locally) the boundary of four minimal surfaces which in general meet at various angles.
Condition (i) ensures that one can prescribe small angles between the opposing minimal surfaces,
smoothly varying along $\Cunder$,
and then realize these angles by appropriately perturbing $\Cunder$ as before.
Since these angles determine the resultant transverse force
(or rather the transverse force density along $\Cunder$)
exerted by $\What$ to $\Cunder$,
the name ``unbalancing condition'' is justified in analogy with other gluing constructions
(see \cite{kapouleas:survey} for a detailed discussion).

Note that we have phrased the theorem in the case that we require the constructed minimal surface
to have the same boundary as the given one, and this is why we impose Dirichlet conditions on $\partial\W$
in (i) and (ii).
We could require some other boundary conditions as well,
and then (i) and (ii) would be appropriately modified.
We also remark that if a symmetry group $\group$ is imposed on the construction,
then conditions (i-ii) should be interpreted to apply only to functions invariant under the action of
$\group$.
This is often equivalent to restricting to functions invariant under the action of $\groupsystem$
which is a larger group and independent of the number of handles used in the desingularization.
An example where the construction fails because these conditions fail,
as we discuss in Section \ref{Ssph},
is provided by two equatorial two-spheres intersecting non-orthogonally in the unit three-sphere.

\subsection*{Open questions and further constructions}
$\phantom{ab}$
\nopagebreak

We briefly comment now on possible extensions of \ref{Tmain}.
Such extensions are related to the possible weakening of conditions (a-d):
Condition (a) can be removed on a case-by-case basis and
actually already in \cite{kapouleas:finite} does not apply.
Having non-compact $\Cunder$, as in the case of two intersecting catenoids
whose axes are parallel and their waists coplanar but not intersecting,
requires a careful understanding of the behavior when the angle of the Scherk surfaces
degenerates to $0$,
and seems possible with further work.
Removing condition (d) seems possible with further work and seems necessary for applications
to the embedded case of the Calabi-Yau problem for minimal surfaces.

Removing condition (c) also involves the degeneration of the angle to $0$ and seems possible with further work.
There are interesting applications of constructions where (c) does not apply,
as for example in desingularizing the intersection of a Clifford torus and an equatorial two-sphere in the round three-sphere.
Finally removing condition (b) seems the hardest, at least in the case when there is a triple intersection on a single component
of $\Cunder$. (Otherwise one can attempt to desingularize on the components of $\Cunder$ successively).
This seems to require a new model for the handles used.
For such a model 
one would try to find and carefully study
minimal surfaces desingularizing three intersecting planes and asymptotic to singly periodic Scherk surfaces
close to the lines of intersection at $\infty$.

\section{Minimal surfaces in the round three-sphere}
\label{Ssph}
\nopagebreak

\subsection*{Constructions of closed embedded minimal surfaces in $\Sph^3(1)$}
$\phantom{ab}$
\nopagebreak

We discuss now known and potential constructions of closed embedded minimal surfaces in $\Sph^3(1)$.
The simplest such construction would involve desingularizing two equatorial two-spheres
intersecting along a great circle.
The angle between the two spheres is the only free parameter of this system.
Unfortunately both conditions (i-ii) in \ref{Tmain} fail:
(i) is violated by the first harmonics on 
each of the intersecting spheres;
(ii) is violated by the first harmonic vanishing on the boundary on each of the four
hemispheres constituting $\What$.
This could be corrected by imposing enough symmetry.
However when the equatorial spheres are not orthogonal,
the first harmonics vanishing on the intersection circle survive all symmetry
and still violate conditions (i-ii).

When the equatorial spheres intersect orthogonally,
$\groupsystem$ is much larger because it includes the reflections through
great circles contained on the two-spheres and orthogonal to the intersection circle.
A group of symmetries $\group$ can be imposed then on the construction which includes
some of these reflections.
These reflections correspond to the extra symmetries the Scherk surfaces possess when their asymptotic
half-planes are orthogonal, and which are reflections through straight lines (see \ref{PScherk}.i).
The first harmonics which previously obstructed conditions (i-ii),
are no longer allowed as a result of
these extra symmetries;
therefore Theorem \ref{Tmain} applies to our system
when the two equatorial spheres intersect orthogonally.
In applying \ref{Tmain} we are using only a very special (and much easier) case of the general theorem.
This follows as we not only 
have a highly symmetric construction (see section \ref{Sdesingularization}),
but also we have so much symmetry that the extended kernel is trivial.
In fact many of the difficulties of even the highly symmetric case become irrelevant.

The surfaces obtained have the same symmetries,
and one expects that they are exactly the same,
as the Lawson surfaces $\xi_{m,k}$ \cite {L2} for $m$ large and $k=1$.
A similar gluing construction should produce surfaces with the same symmetries,
and one expects identical,
to the Lawson surfaces $\xi_{m,k}$ for $m$ large and $k>1$.
Such a construction would be considered a desingularization of $k+1$ spheres symmetrically
arranged around a common great circle.
The Karcher-Scherk towers \cite{Ka1} which desingularize $k+1$ planes
that are symmetrically arranged around a common line of intersection
provide an appropriate model to be used instead of the Scherk surfaces.

By incorporating Clifford tori and rotationally invariant Delaunay-like surfaces
one can propose various new desingularization constructions.
Most of them are highly symmetric,
but there are also some which are much more demanding,
as for example the desingularization of a Clifford torus intersecting a great two-sphere,
where the intersection curve contains points where the torus and the sphere are tangent to each other.
We intend to study these constructions elsewhere.
Desingularizations of intersecting Clifford tori can be constructed also by the original approach of Lawson
(work in progress of Choe and Soret \cite{jaiyoung}).

We have already discussed in Section \ref{Sdoubling} doubling constructions for the equatorial two-sphere
and the Clifford torus.
We discuss now a ``second generation'' construction where doublings of Clifford tori are combined and desingularized
to produce new closed embedded minimal surfaces.
We start by assuming given a finite sequence of natural numbers $k_1< ... <k_r$,
where $k_1=1$,
and if $i<j$ then $k_j$ is a multiple of $k_i$.
We assume then $m$ is large enough and
we consider the square lattices 
$$
L_1\subset L_2\subset ... \subset L_r
$$
where $L_i$ is an $m_i\times m_i$ lattice,
where $m_i:=k_i m$.
Since $m$ is large enough,
\ref{Tclifford} applies,
and we obtain doublings $\Tbb_i$ of configuration $L_i$ ($i=1,...,r$).

We take $\W$ to be the collection of the $\Tbb_i$'s,
and we may or may not include also the Clifford torus $\Sigma$ itself.
We take $\Cunder$ to consist of the intersection points.
By \ref{Etb} the ratios of the sizes of the catenoidal bridges are either
very small or very large.
It is easy to check then that $\Cunder$
is the disjoint union of smooth circles of transverse intersection.
We want to desingularize then $\W$ along $\Cunder$,
where we impose as symmetry group $\group$ for the desingularization 
the group of isometries of $L_1$ which is also the group of symmetries for the construction
of the doubling $\Tbb_1$, but strictly smaller than the group of symmetries for $\Tbb_i$ when $i>1$.

\addtocounter{equation}{1}
\begin{theorem}[\cite{kapouleas:nested}]
\label{Tnested}
If $m$ is large enough depending on the sequence $\{k_i\}_{i=1}^r$,
then $\W$ constructed as above can be desingularized along $\Cunder$
with imposed group of symmetries $\group$
as in Theorem \ref{Tmain},
because all conditions (a-d) and (i-ii) apply.
We obtain then a family of closed embedded minimal connected surfaces
parametrized by the number of handles used in desingularizing
each component of $\Cunder$ subject to the action of $\group$.
\end{theorem}

Note that the picture in this theorem is similar to some extent to the one for the
desingularization theorem in \cite{kapouleas:finite}:
In the vicinity of a point $p\in L_i$ we have approximately coaxial catenoidal bridges 
and (if $\Sigma$ is included in $\W$) an approximate plane through the waist.
The rotational symmetry is only approximate however, and completely fails further away from $p$.
Actually one can check that in general there are components of $\Cunder$
which have trivial stabilizer in $\group$.
The full force of \ref{Tmain} is therefore needed.

To prove Theorem \ref{Tnested} we have to check conditions (a-d) and (i-ii) apply.
(d) is automatic because there is no boundary,
and (a-c) follow easily because the geometry of the $\Tbb_i$'s
is well controlled and the ratios of the sizes of the catenoidal bridges are extreme.
Proving conditions (i-ii) is harder and requires the use of the ``geometric principle''.
Condition (i) actually
amounts to non-existence of 
zero eigenvalues for the linearized
operator $\Lop$ on a doubling of the Clifford torus constructed as in Theorem \ref{Tclifford}
not only when we restrict to functions invariant under the action of the symmetry group $\group_{doubling}$
which was imposed on the construction,
but also when we restrict to functions 
invariant under the strictly smaller
symmetry group $\group'_{doubling}$ of a strictly smaller square sublattice of the lattice used in the doubling
construction.
This is because in the desingularization construction we only impose the group $\group$ which
is the symmetry group of the lattice $L_1$,
while on the doubling construction for $\Sigma_i$ we impose the larger symmetry group of the larger lattice $L_i$.

To facilitate the discussion we assume that $\Tbb$ is the doubling of the Clifford torus as in \ref{Tclifford}
with configuration an $m'\times m'$ square lattice and we consider functions invariant under the symmetries of a smaller
$m \times m$ square sublattice.
The fundamental domain by the rotational symmetries has now $(m')^2/m^2$ catenoidal bridges.
By conformality the kernel under the new symmetries is nontrivial
if and only if the kernel for $\Lh$ (defined as in \ref{ELh}) is.
Since many of the catenoidal bridges remain free by the symmetries to ``slide'' to new positions
the approximate kernel has some finite dimension which is a function of $m'/m$.
We apply now the ``geometric principle'':
By imposing appropriate ``slidings'' or ``relocations'' on the standard regions,
and transiting to solutions to the linearized equation on the transition regions,
with the corresponding Dirichlet data on their boundaries,
we can obtain functions on the surface which can be analyzed to demonstrate
by integrating by parts that the small eigenvalues never vanish.

To establish that condition (ii) also holds,
consider the domains into which $\W$ is subdivided by $\Cunder$.
Because of the extreme (either very small or very large) ratios of the sizes of the catenoidal bridges,
these domains are either very narrow, or are large enough to be analyzed
in a similar way as the surfaces were to establish condition (i).
In either case we avoid zero eigenvalues and the result can be proven.

\subsection*{Can the Lawson surfaces flap their wings and other open questions}
$\phantom{ab}$
\nopagebreak

A general question asks for the classification of all closed minimal surfaces in $\Sph^3(1)$
whose area is less than $4\pi C$, where $C$ is some given constant.
Even if $C$ is a small integer, 
this is a very hard question since it includes the following two questions:
The Lawson conjecture that the Clifford torus is the only embedded minimal torus and
the question (recommended to the author by Mark Haskins)
whether the Clifford torus has the smallest area
among closed embedded minimal surfaces with the exception of the equatorial two-sphere.
We do not discuss these questions further but indicatively we mention \cite{ros:lc,ros:area,white}.

The only closed minimal surfaces in $\Sph^3(1)$ known to the author,
besides the equatorial sphere, the Clifford torus, and surfaces constructed by gluing methods,
are the Lawson surfaces ($\xi_{m,k}$ in the notation of \cite{L2}),
and the surfaces constructed in \cite{KPS}.
As we have discussed already these surfaces can be considered as desingularizations or doublings,
and so are the surfaces discussed so far by gluing methods.
So it would be of some interest to answer the following,
where by ``non-elementary'' we mean not the equatorial sphere or the Clifford torus.

\addtocounter{equation}{1}
\begin{question}
\label{Qdifferent}
Find a non-elementary connected closed minimal surface in $\Sph^3(1)$ which is not a doubling or a desingularization.
\end{question}

Another question which seems important is to determine 
to what extent the Lawson surfaces are unique,
or more precisely characterized by some simple properties they possess.
We formulate such a question in \ref{QLone} below for the case $k=1$ in the notation of \cite{L2}.
We first give the following definitions:
We define the $\delta$-neighborhood $X_\delta$ of a set $X\subset\Sph^3(1)$
to be the set of points whose distance from $X$ is $\le\delta$.
We also define a $\delta$-wing to be a minimal disc in $\Sph^3(1)$
which can be considered as the graph of a function $f$ over a domain of an equatorial two-sphere,
and where $\|f:C^0\|\le\delta$.

\addtocounter{equation}{1}
\begin{question}
\label{QLone}
Is there a $\delta>0$ such that the following holds?
If a closed embedded minimal surface $M$ in $\Sph^3(1)$ has the property that there is a great circle $\C$,
such that $M\setminus \C_\delta$ consists of four $\delta$-wings,
then $M$ is one of the Lawson surfaces $\xi_{m,k}$ with $k=1$.
\end{question}

This question seems related to uniqueness questions for singly periodic Scherk surfaces
in which there has been some recent progress 
\cite{meeks:unique}.
An affirmative answer to Question \ref{QLone} would imply 
in particular that high-genus Lawson surfaces cannot ``flap'' their wings.
It would also imply that two intersecting non-orthogonal equatorial spheres
cannot be the varifold limit of a sequence of desingularizations.
Similarly an equatorial sphere covered twice cannot be the limit of a sequence
of doublings where all the catenoidal bridges concentrate close to an equatorial circle.

The following heuristic argument provides evidence for the last two statements,
and hence for an affirmative answer to \ref{QLone} as well.
We concentrate on the case of the intersecting spheres since the argument in the other case is similar.
Let $\C$ be the circle of intersection of the two spheres
and $M$ a desingularization close to the limit.
It is reasonable to assume that for some small $\delta$,
$M\cap \C_\delta$ is approximated very closely by a Scherk surface bent along $\C$
and scaled to a size $\tau$ much smaller than $\delta$.
This motivates us to assume also that $M\setminus \C_{c\tau}$,
where $c$ is a large constant but such that $c\tau$ is much smaller than $\delta$,
consists of four discs each of which can be described as a graph over domains of the given hemispheres.
More precisely one of these discs is a graph of a function $f$ over a rotationally invariant domain $\Omega$,
such that $\partial\Omega$ is a parallel circle to $\C$.
Because of the smallness of $f$ we can assume that it approximately satisfies the linearized equation
$$
\Delta f + 2 f=0,
$$
and moreover that $f-\overline{f}$ is much smaller than $f$ and decays fast away from $\partial\Omega$,
where $\overline{f}$ is defined by
$$
\left.\overline{f}\right|_{\C'}
=
\avg_{\C'} f
$$
on each circle $\C'$ parallel to $\C$.

We have then that $\overline{f}$ also satisfies the linearized equation above,
which now amounts to an ODE.
By uniqueness of ODE solutions $\overline{f}$ is a multiple of the first harmonic on the sphere containing $\Omega$
and vanishing on $\C$.
This implies that $|\nabla f|$ is close to $|f|/\delta$ on $\Omega\cap\partial\C_\delta$.
We have assumed however that in $\C_\delta$ the surface approximates closely a Scherk surface
and therefore by \ref{PScherk}.ii on $\Omega\cap\partial\C_\delta$
$f$ is close to $b_\alpha\tau$ and by the exponential decay $|\nabla f|$ is of order $e^{-\delta/\tau}$.
By assuming also that $\tau$ is small enough in terms of $b_\alpha$,
which we expect for the surfaces in the sequence that are close enough to the intersecting spheres,
we conclude that unless $b_{\alpha}$ vanishes,
$|\nabla f|$ is much smaller than $|f|$ on $\Omega\cap\partial\C_\delta$.
This contradicts the previous conclusion that 
$|\nabla f|$ is close to $|f|/\delta$ on $\Omega\cap\partial\C_\delta$.
Thus $b_{\alpha}=0$,
which amounts to the spheres being orthogonal.

Morally the argument above depends upon the following.
In Euclidean space we have a two-parameter family of planes parallel to the axis of the Scherk,
and the asymptotic half-planes of the Scherk surfaces do not contain the axis when they are not orthogonal.
In the three-sphere once we determine the circle of intersection $\C$,
we have only a one-parameter family of spheres ``parallel'' to $\C$.
These spheres always contain $\C$.
This creates an incompatibility unless the limit spheres are orthogonal.

Assuming now that the answer to Question \ref{QLone} is affirmative,
there are two different directions we can extend the question.
First, we can ask if we start increasing $\delta$ for what value and how the statement will fail.
Second, we can extend the question by considering the case where $M\setminus \C_\delta$ consists of $d$ discs instead of four.
In \ref{QLtwo} we formulate one of the stronger such questions that is possible to ask.
(Note that the assumptions on the components of $M\setminus \C_\delta$ are weaker than in \ref{QLone}.)
If the answer to \ref{QLtwo} or a modification of it is positive,
one could ask furthermore how the ``best'' $\delta$ depends on $d$.
In less formal language \ref{QLtwo} asks whether the Lawson surfaces are characterized by the fact
that their topology concentrates near a great circle:

\addtocounter{equation}{1}
\begin{question}
\label{QLtwo}
Is there a $\delta>0$ such that the following holds?
If a closed embedded minimal surface $M$ in $\Sph^3(1)$ has the property that there is a great circle $\C$,
such that $M\setminus \C_\delta$ is the union of minimal discs,
then $M$ is one of the Lawson surfaces $\xi_{m,k}$ described in \cite{L2}.
\end{question}

Similar questions to the above can be asked for the desingularizations of intersecting Clifford tori.

\section{The building blocks for the desingularization construction}
\label{Sbuildingblocks}
\nopagebreak

\subsection*{The given system of minimal surfaces}
$\phantom{ab}$
\nopagebreak

As we discussed earlier we assume given minimal surfaces intersecting transversely along
a simple closed curve along which we intend to desingularize.
We describe this now in a precise and systematic way.
We can think of the intersecting minimal surfaces as the components of a single
minimal surface.
More precisely we assume that $\W$ is an abstract surface which is immersed into a
Riemannian manifold $\N$ by a minimal immersion $\X:\W\to\N$.
$\W$ may or may not have boundary, and we denote the metric of the Riemannian manifold by $g$.
$\C$ is a smooth closed curve embedded in the interior of $\W$,
and the restriction of $\X$ to $\C$ is a double covering of a curve $\Cunder$ in $\N$.
$\Cunder$ is the curve of intersection to be desingularized.

\addtocounter{equation}{1}
\begin{convention}
\label{CCunder}
To simplify the presentation we will treat $\Cunder$ as an embedded curve in $\N$.
Obvious modifications however would allow an immersed $\Cunder$.
\end{convention}

Locally we can view $\Cunder$ as the boundary of four minimal pieces (described by $\X$),
which can be perturbed independently so that they have new boundaries.
To describe this precisely we need to define a new compact abstract surface with boundary $\What$,
which is obtained from $\W$ by ``cutting'' along $\C$.
We denote by $\QQ:\What\to\W$ the smooth map which is the ``natural projection'' of $\What$ to $\W$:

\addtocounter{equation}{1}
\begin{definition}
\label{DWhat}
$\What$ and $\QQ:\What\to\W$ are characterized by the following:
\newline
(i). $\Chat:=\QQ^{-1}(\C)\subset\partial\What$.
\newline
(ii). $\left.\QQ\right|_{\What\setminus\Chat}$ is a diffeomorphism
onto $\W\setminus\C$.
\newline
(iii). The restriction of $\QQ$ to $\Chat$,
$\left.\QQ\right|_\Chat$,
gives a two-to-one covering of $\C$ by $\Chat$.
\newline
(iv).
$\QQ$ is a quotient map.

We also define $\Xhat:=\X \circ \QQ:\What\to\N$.
\end{definition}

Note that $\partial\What$ is the disjoint union of $\Chat$ and $\QQ^{-1}(\partial\W)$.
Note also that $\Xhat$ is a minimal immersion.
$\Xhat$ restricted to $\Chat$ provides a four-to-one covering of $\Cunder$.
At each point $p\in\Cunder$ then, $\Xhat$ provides four (unit inward)
conormals, one at each of the
four points of $\Chat$ mapping to $p$.
We denote the corresponding vectors in $T_p\N$ by
\addtocounter{theorem}{1}
\begin{equation}
\label{Evectors}
\vece_{1}(p),\quad
\vece_{-1}(p):=-\vece_1(p),\quad
\vece_{2}(p),\quad
\vece_{-2}(p):=-\vece_2(p),
\end{equation}
and we call them the conormal vectors to $\Cunder$ at $p$.
Note that the choice (among the four vectors)
of $\vece_1$ is arbitrary, 
and then the choice (among the remaining two vectors) of $\vece_2$ is arbitrary.
We remark that later we will perturb $\Xhat$ so that the four conormal vectors are not opposite
in pairs anymore.

Since $\Cunder$ is a closed curve, we have 
\addtocounter{theorem}{1}
\begin{equation}
\label{ECj}
\Cunder=\bigcup_{j=1}^k\Cunder_j,
\end{equation}
where $k$ is the (finite) number of connected components,
and each $\Cunder_j$ a smooth circle.
We will need precise notation to deal with tubular neighborhoods of these
circles, and the conormals of the minimal pieces meeting along them.
For each of these circles $\Cunder_j$ we define
in \ref{DAj} below
a smooth map
$\A_j:\Rbb^3\to N\Cunder_j$,
where $N\Cunder_j$ is the normal bundle of $\Cunder_j$,
and compatible unit-speed parametrizations
$\al_j:\Rbb\to\Cunder_j$ 
and
$\alhat_{j,i}:\Rbb\to\Chat$.
It is useful to fix some helpful notation first:

\addtocounter{equation}{1}
\begin{notation}
\label{NR}
We use $R_1$ to denote the identity map on $\Rtwo$,
and $R_2$, $R_{-1}$, and $R_{-2}$ the reflections in $\Rtwo$
with respect to the $x_1$-axis, the origin, and the $x_2$-axis respectively.
We also let each $R_i$ act on $\Rthree$ by $R_i(x_1,x_2,x_3)=(R_i(x_1,x_2),x_3)$.
\end{notation}

\addtocounter{equation}{1}
\begin{definition}
\label{DAj}
$\A_j:\Rbb^3\to N\Cunder_j$,
$\al_j:\Rbb\to\Cunder_j$,
and
$\alhat_{j,i}:\Rbb\to\Chat$,
($\iequal$)
are defined to be smooth maps 
characterized by the following:
\newline
(i).
$\A_j$ is a bundle isometry if we consider $\Rbb^3$ as a trivial bundle over the $x_3$-axis equipped with the
standard Euclidean metric.
\newline
(ii).
The restriction of $\A_j$ to the $x_3$-axis is a unit-speed covering of $\Cunder_j$
which we denote by $\al_j:\Rbb\to\Cunder\subset\N$.
We denote by $\ell_j$ the length of $\Cunder_j$ and we have then
\addtocounter{theorem}{1}
\begin{equation}
\label{Elk}
\al_j(x+\ell_j)  \equiv   \al_j(x).
\end{equation}
(iii).
In the fiber above each point $(0,0,x)$ there are four vectors
(recall \ref{NR}),
$$
R_i(\cos\alpha_j(x),\,\sin\alpha_j(x),\,x)
\qquad
(\iequal),
$$
which map by $\A_j$ to the four unit conormal vectors at $\al_j(x)$.
Note that the angle of intersection at $\al_j(x)$ is then $2\alpha_j(x)$ or $\pi-2\alpha_j(x)$.
\newline
(iv).
For $x\in\Rbb$ and $\iequal$ we have that
$$
\Xhat_*\circ\veceta_\ji(x)
=
\A_j(  R_i(\cos\alpha_j(x),\,\sin\alpha_j(x),\,x)  ),
\qquad
\text{where}
\qquad
\veceta_\ji(x):=
\left.\veceta\right|_{\alhat_\ji(x)} ,
$$
where $\left.\veceta\right|_p \in T_p\What$ denotes the inward unit normal to $\Chat$ at some $p\in\Chat$.
\newline
(v).
The above characterize $\A_j$ up to arbitrary choices of the following which we assume from now on:
First, an arbitrary point on
$\Cunder_j$ to serve as $\al_j(0)$.
Second, a choice of the orientation of $\Cunder_j$ induced by $\al_j$.
Third, 
a choice of $\alhat_{j,1}(0)$ among the four possible points,
and then a choice of $\alhat_{j,2}(0)$ among the remaining two.
\end{definition}

Note that because we assume the intersection to be transverse we have
$\alpha_j\in(0,\pi/2)$.
By compactness this implies that 
\addtocounter{theorem}{1}
\begin{equation}
\label{Edel0}
\alpha_j\in(20\delta_0,\frac\pi2-20\delta_0)
\end{equation}
for some $\delta_0>0$ which depends on the given system
of minimal surfaces.
Note also that $\A_j$ satisfies
\addtocounter{theorem}{1}
\begin{equation}
\label{EAlk}
\A_j(x_1,x_2,x+\ell_j)  \equiv   \A_j(P_j(x_1,x_2), x),
\end{equation}
where $P_j$ is a Euclidean motion (or the identity map) on the $x_1x_2$-plane
which fixes the union of the coordinate axes.
An example where $P_j$ is not the identity is provided by a minimal Moebius band in Euclidean three-space
with a circle of self-intersection.

\subsection*{The singly periodic Scherk surfaces}
$\phantom{ab}$
\nopagebreak

We discuss now the singly periodic Scherk surfaces (recall \ref{EScherk}),
which we use as ingredients in our construction.
Because $\Sa$ degenerates as $\alpha\to0,\pi/2$,
and we need uniform bounds on the geometry of the $\Sa$'s,
we restrict $\alpha$ away from $0$ and $\pi/2$:
We assume from now on that in accordance with \ref{Edel0}
\addtocounter{theorem}{1}
\begin{equation}
\label{Ealpha}
\alpha\in[10\delta_0,\frac\pi2-10\delta_0].
\end{equation}
Before we discuss further the Scherk surfaces, we adopt some helpful notation:
\addtocounter{equation}{1}
\begin{definition}
\label{Dvece}
We define $\vece_{i,\alpha}$ and $\veceperp_{i,\alpha}$ for $\iequal$
and $\alpha$ as in \ref{Ealpha} by (recall \ref{NR})
$$
\vece_{1,\alpha}:=(\cos\alpha,\sin\alpha,0),
\quad
\veceperp_{1,\alpha}:=(-\sin\alpha,\cos\alpha,0),
\quad
\vece_\ialpha:=R_i(\vece_\onealpha),
\quad
\veceperp_\ialpha:=R_i(\veceperp_\onealpha).
$$
\end{definition}

\addtocounter{equation}{1}
\begin{definition}
\label{DgroupS}
We define $\groupS$ to be the group of Euclidean motions in $\Rbb^3$
generated by the reflections with respect to the coordinate planes and the plane
$\{x_3=\pi\}$.
\end{definition}

In the next proposition we enumerate the properties of the
Scherk surfaces which are relevant to our constructions.
Note that these surfaces can be decomposed into a ``core'',
which is within a finite distance from the $x_3$-axis,
and four ``wings''
which are asymptotic to four half-planes
symmetrically arranged around the $x_3$-axis.
The wings depend up to Euclidean motion only on the parameter $\alpha$
of the Scherk surface in consideration.
  
\addtocounter{equation}{1}
\begin{prop}
\label{PScherk}
$\Sa$ is a singly periodic
embedded complete minimal surface which
depends smoothly on $\alpha$ and has the following properties:
\newline
(i). $\Sa$ is invariant under the action of $\groupS$.
$\Sc:=\Scortho$ 
has a larger symmetry group containing 
also the reflections with respect to the lines
$\{x_1=\pm x_2,x_3=\pi/2\}$.
\newline
(ii).
If we exclude a neighborhood of the $x_3$-axis the rest of $\Sa$
has four components we call the ``wings'' of $\Sa$.
Each wing can be described as a graph $\BB_\ialpha(\Rtwoplus)$
over an asymptotic half-plane $\bb_\ialpha(\Rtwoplus)$,
where $\Rtwoplus:=\{(x,s):s\ge0\}$, 
$\iequal$,
and
$\bb_\ialpha,\BB_\ialpha:\Rtwoplus\to\Rthree$
are defined by
$$
\BB_\ialpha(x,s)=\bb_\ialpha(x,s)+\fa(x,s)\,\veceperp_{i,\alpha},
\qquad
\bb_{i,\alpha}(x,s)=(0,0,x)+(s_0+s)\vece_{1,\alpha}+b_\alpha \veceperp_{i,\alpha},
$$
where 
$b_\alpha:=\sin2\alpha\,\log(\cot\alpha)$,
$s_0>0$ is a constant which depends only on $\delta_0$ and a given 
$\varepsilon\in(0,10^{-3})$,
and $\fa:\Rbb^2_+\to[-\varepsilon,\varepsilon]$
is a smooth function
which depends only on $\alpha$.
Moreover
$\fa$, $\bb_{i,\alpha}$, and $\BB_\ialpha$,
depend smoothly on
$\alpha$ as in \ref{Ealpha},
and (iii)-(vi) below are satisfied.
\newline
(iii).
$\fa$ is periodic in the sense that $\fa(x,s)\equiv\fa(x+2\pi,s)$.
\newline
(iv).
$\Sa\setminus\cup_{i=1}^4\BB_\ialpha(\Rbb^2_+)$,
which we call the core of $\Sa$,
is connected and lies within distance $2s_0$ from the $x_3$-axis.
\newline
(v).
$\|\fa:C^k(\Rbb^2_+,e^{-s}) \|\le C(k)\,\varepsilon$
and
$\|d\fa/d\alpha:C^k(\Rbb^2_+,e^{-s}) \|\le C(k)\,\varepsilon$.
\newline
(vi). $ |b_\alpha| + |db_\alpha/d\alpha| < \varepsilon s_0$.
(Notice that the right hand side is not small because $s_0$
is large.)
\end{prop}

It is important to understand the Gauss map $\gm:\Sa\to\Stwo$
and the geometry of the pullback metric
$h:=\gm^*g_{\Stwo}$.
We can decompose $\Sa$ into a sequence of regions
$$
S_n=\Sa\cap\{x_3\in[n\pi,(n+1)\pi]\}
$$
which are permuted transitively by the symmetry group $\groupS$ of $\Sa$.
It is enough then to study one of them as in the next proposition.
It follows that $(\Sa,h)$ is
an isometric cover of
$\Stwo\setminus\{(\veceperp_\ialpha\}$
with covering map the Gauss map $\gm$.
Each $S_n$ can be identified by $\gm$
with a closed hemisphere with four points
on its boundary circle removed.
Because of this we call the $S_n$'s the hemispherical regions of $\Sa$.

Each $\partial S_n$ under this identification is the disjoint
union of four open arcs.
$S_n$ connects to the adjacent hemispheres $S_{n+1}$ and $S_{n-1}$
through these arcs so that
$S_n\cap S_{n+1}$ consists of two opposing arcs 
and
$S_n\cap S_{n-1}$ consists of the other two.
Each of the four points removed corresponds to the $\infty$
of one of the four wings.
Each wing of $\Sa$ is the preimage under the Gauss map $\gm$
of a neighborhood in $\Stwo$ of such a point.

\addtocounter{equation}{1}
\begin{prop}
\label{Pgauss}
The Gauss map $\gm$ of $\Sa$ has the following properties:
\newline
(i).
$\gm$ restricted to
$S_0=\Sa\cap\{x_3\in[0,\pi]\}$
is a diffeomorphism onto
the closed hemisphere with four boundary points removed
$\Stwo\cap\{x_3\ge0\}\setminus\{\veceperp_\ialpha\}$.
\newline
(ii). Let $E_i$ ($i=1,...,4$)
be the arcs into which the equator
$\Stwo\cap\{x_3=0\}$ is decomposed by removing the points 
$\{\veceperp_\ialpha\}$,
numbered counterclockwise so that $(1,0,0)\in E_1$.
We have then
$$
\gm(\Sa\cap\{x_3=0\})=E_1\cup E_3,
\qquad
\gm(\Sa\cap\{x_3=\pi\})=E_2\cup E_4.
$$
(iii). $\Sa$ has no umbilics and $\gm^*g_{\Stwo}=\frac12|A|^2g$.
\end{prop}

In order to describe the constructions later we need appropriate parametrizations of the
Scherk surfaces which we describe in the next definition.
Note that
by condition (iii) below,
the new parametrizations $Z_\alpha$ are consistent with the parametrizations
of the wings 
$\BB_\ialpha$
defined in \ref{PScherk}.

\addtocounter{equation}{1}
\begin{definition}
\label{DZalpha}
We fix smooth embeddings
$ Z_\alpha:\Sc\to\Rthree$
($\Sc:=\Scortho$ 
and 
$\alpha$ is as in \ref{Ealpha})
satisfying the following:
\newline
(i).
$ Z_\alpha$ depends smoothly on $\alpha$ and
$Z_\alpha(\Sc)=\Sa$.
\newline
(ii).
$ Z_\alpha$ is equivariant under the action of $\groupS$.
\newline
(iii).
$ Z_\alpha\circ
\BB_{i,\pi/4}=
\BB_\ialpha$.
\newline
(iv).
There is a $\delta_1>0$, which depends only on $\delta_0$,
such that if $p\in\Sc$ and
$|x_3(p)-n\pi|<\delta_1$
for some $n\in \mathbb{Z}$,
then
$x_3\circ Z_\alpha(p)=x_3(p)$.
\end{definition}

\subsection*{Dislocations of the Scherk surfaces}
$\phantom{ab}$
\nopagebreak

In the construction of our initial surfaces we need to allow for certain
dislocations as required by the general methodology.
These dislocations force perturbations
of the Scherk surfaces which introduce mean curvature 
in accordance with the
``geometric principle''.
The modifications of the Scherk surfaces are controlled by seven parameters
as follows:
$\sigma\in\Rbb$ controls the rate of change of scale which is related to the
creation of ``longitudinal kernel'',
that is the kernel induced by translations in the direction of the $x_3$-axis.
$\vartheta:=(\vartheta_1,\vartheta_2)\in\Rtwo$ controls dislocations where opposing
wings are rotated relative to each other in order to create ``transverse kernel'',
that is kernel induced by translations perpendicular to the $x_3$-axis.
Finally 
$\varphi:=(\varphi_i)_{\iequal}\in\Rfour$
introduces relative rotations between each wing and the core.
This corresponds to the creation of extended substitute kernel
required for arranging the decay of solutions to the linearized equation
along the wings.

Since $\sigma$ relates to the rate of change of scale,
we postpone discussing it until the next section,
where we study modifications of the Scherk surfaces
where the controlling parameters vary along the surface.
In this subsection then we will only study how to
modify the embedding $Z_\alpha$ of $\Sa$
to an embedding $Z_\althph$ for the modified surface.
We assume that $\alpha$ is as in \ref{Ealpha},
and $\vartheta$ and $\varphi$ satisfy
\addtocounter{theorem}{1}
\begin{equation}
\label{Ethph}
|\vartheta|\le\delta_0, \qquad |\varphi|\le\delta_0.
\end{equation}

\addtocounter{equation}{1}
\begin{definition}
\label{DDvartheta}
We fix a family of diffeomorphisms $D_{\vartheta}:\Rtwo\to\Rtwo$
which have the following properties:
\newline
(i). $D_{\vartheta}$ depends smoothly on $\vartheta$.
\newline
(ii). For $\vartheta=(0,0)$ $D_{\vartheta}$ is the identity.
\newline
(iii).
$D_\vartheta$ restricted to a neighborhood of the coordinate lines is the identity.
\newline
(iv).
On the region
$R_i\{(r\cos\phi,r\sin\phi): r>5,\,\,
\phi\in\,[\pi+5\delta_0,\,\frac{3\pi}2- 5\,\delta_0]\}
$
(recall \ref{NR}),
$D_\vartheta$ acts as a rotation around the origin
by an angle
$\vartheta_1/2$ if $i=1$,
$-\vartheta_1/2$ if $i=-1$,
$\vartheta_2/2$ if $i=2$,
and
$-\vartheta_2/2$ if $i=-2$.
\newline
(v).
$R_{-1}\circ D_{\vartheta}=D_{-\vartheta}$.
\end{definition}

We extend the action of $D_\vartheta$ to $\Rthree$
by requiring that it leaves the $x_3$ coordinate unchanged.
$D_\vartheta(\Sa)$ is then a modification of $\Sa$ where opposing wings are
rotated in a symmetric way to create relative angles 
$\vartheta_1$ and $\vartheta_2$ respectively.
More precisely we have
\addtocounter{theorem}{1}
\begin{equation}
\label{EDwings}
D_\vartheta\circ Z_\alpha\circ \BB_{i,\pi/4}=\RRR_i\circ\BB_\ialpha,
\end{equation}
where $\RRR_i$ is a rotation by an angle $\pm\vartheta_{|i|}/2$.

\addtocounter{equation}{1}
\begin{definition}
\label{DZalthph}
We define now $Z_\althph:\Sc\to\Rthree$ by requiring the following:
\newline
(i).
On the core of $\Sc$ we have $Z_\althph=D_\vartheta\circ Z_\alpha$
(recall \ref{PScherk}.iv).
\newline
(ii).
On $\{s\ge1\}\subset\Rtwoplus$ we have 
(recall \ref{DZalpha}.iii)
$Z_\althph\circ\BB_{i,\pi/4}=\RRR'_i\circ\RRR_i\circ\BB_\ialpha$,
where $\RRR'_i$ is the rotation around the line $\RRR_i\circ\bb_\ialpha(\{s=0\})$
by an angle $\varphi_i$.
\newline
(iii).
It remains to match (i) where $\RRR'_i$ is not applied,
to (ii) where it is.
To define $Z_\althph$ on $\BB_{i,\pi/4}(\Rbb\times[0,1])$
we require that for $s\in[0,1]$ we have
\begin{multline*}
Z_\althph\circ\BB_{i,\pi/4}(x,s)=
\\
(1-\psi[0,1](s)\, )
D_\vartheta\circ Z_\alpha\circ\BB_{i,\pi/4}(x,s)
+
\psi[0,1](s)\,
\RRR'_i\circ\RRR_i\circ\BB_\ialpha(x,s).
\end{multline*}
\end{definition}

Note that in 
\ref{DZalthph}.iii
we have a smooth transition from the core which has not been rotated by $\RRR'_i$,
to the wing which has been.
We need to introduce one last modification by simply rotating around the $x_3$-axis:
We define
\addtocounter{theorem}{1}
\begin{equation}
\label{EZfinal}
Z_\althphp:=\RRR''\circ Z_\althph,
\end{equation}
where $\RRR''$ is the rotation around the $x_3$-axis by a
given angle $\phi\in [-\delta_0,\delta_0]$.

For future reference it is more convenient to replace the parameters of $Z_\althphp$ with the asymptotics
of the surface and the $\varphi_i$'s as follows.
Let $\RRR'''_i$ denote the rotation in $\Rtwo$ around the origin by an angle $\phi+\varphi_i$.
The unit vector $\vecep:=\RRR'''_i\circ\RRR_i(\vece_\ialpha)$
is parallel to the asymptotic half plane of the $i$-th wing of the modified
Scherk surface.
We call the collection of these unit
vectors and the angles $\varphi_i$ the extended tetrad
\addtocounter{theorem}{1}
\begin{equation}
\label{ET}
T:=\{\vecep,\varphi_i\}_{i=\pm1,\pm2}
\end{equation}
associated to $Z_\althphp$.
It is easy to check that conversely given
$
T=\{\vecepp,\varphi_i\}_{i=\pm1,\pm2}
$
with the angles between $\vecepp$ and $\vecepm$, and the $\varphi_i$'s,
in the interval $[-\delta_0/3,\delta_0/3]$,
there are unique $\althphp$,
such that $T$ satisfies \ref{ET}.
We adopt then the notation
\addtocounter{theorem}{1}
\begin{equation}
\label{EZT}
Z_T:=Z_\althphp.
\end{equation}

\section{An initial surface for the desingularization construction}
\label{Scentral}
\nopagebreak

\subsection*{The construction of the Scherk cores}
$\phantom{ab}$
\nopagebreak

The initial surfaces in the construction for Theorem \ref{Tmain}
form a family parametrized by parameters $\xibo$ where $\xibo$
is an appropriately small element of a finite dimensional vector space of high dimension.
In this section we outline the construction of the ``central'' initial surface,
that is the initial surface with all the parameters $\xibo$ vanishing.
The construction of the Scherk cores is based on
matching 
the cores of the Scherk surfaces
to the geometry of the given minimal surfaces at the vicinity of $\Cunder$
by varying their parameter $\alpha$ along the curve and 
appropriately scaling, bending, and twisting them.
We start by discussing the construction
of the Scherk cores desingularizing $\Cunder$.
We describe this by taking a neighborhood of a Scherk core in the initial surface to be
$\Zti_j(\Score)$,
where
\addtocounter{theorem}{1}
\begin{equation}
\label{EScore}
\Score:=
\Scortho\setminus\cup_{i=1}^4\BB_{i,\pi/4}(\Rbb\times[1/10,\infty)),
\end{equation}
and $\Zti_j:\Score\to\N$ is defined by
\addtocounter{theorem}{1}
\begin{equation}
\label{EEjzero}
\Zti_j:=\E_j\circ\Z_j,
\qquad
\text{where}
\qquad
\E_j:=\exp\circ\A_j\circ\CC_j,
\end{equation}
where $\exp$ is the exponential map of $\N$,
$\A_j$ was defined in \ref{DAj},
$\CC_j:\Rthree\to\Rthree$ is a bundle isomorphism (recall \ref{DAj}.i)
used to adjust the scale,
and $\Z_j:\Sc\to\Rthree$'s construction is based on the embeddings of the
Scherk surfaces discussed in the previous section.
The simplest choices are given by
$$
\CC_j(x_1,x_2,x_3):=\rho_j(x_1,x_2,x_3),
\qquad
\Z_j(p):=Z_{\alpha_j(\rho_j x_3(p))}(p),
$$
where $\rho_j$ is a constant, $p$ is an arbitrary point of $\Sc$,
and $x_3(p)$ is the $x_3$-coordinate of $p\in\Sc\subset\Rthree$.

Note that we cannot avoid varying the $\alpha$ parameter
along $\Sc$ since the core has to match reasonably well the given
minimal surfaces which have a varying angle of intersection along $\Cunder_j$.
Note also that by \ref{Elk} we need that
\addtocounter{theorem}{1}
\begin{equation}
\label{EEjrho}
\Zti_j(\Sc\cap\{x_3=x\})=
\Zti_j(\Sc\cap\{x_3=x+\rho_j^{-1}\ell_j \}).
\end{equation}
This implies in particular that
\addtocounter{theorem}{1}
\begin{equation}
\label{Erhoj}
\rho_j=\frac{\ell_j}{m_j\pi},
\end{equation}
where $m_j$ is a large integer which prescribes the number of hemispherical regions
(equivalently half-handles) that we use to desingularize $\Cunder_j$.
$m_j$ is restricted by the topology of the intersection to be even or odd
but otherwise can be prescribed arbitrarily.
\addtocounter{equation}{1}
\begin{convention}
\label{Cmj}
We assume from now on that each $m_j$ has been chosen,
and it is as large as needed for the construction to work.
To simplify technical aspects of the construction we also assume given a constant $\cunder_0>0$,
and that all the ratios $m_j/m_{j'}$ are bounded by $\cunder_0$.
\end{convention}

Before we proceed we augment the definitions in \ref{DAj} 
by defining smooth maps $\vecnu_\ji:\Rbb\to T\N$ as follows (recall \ref{Dvece}).
\addtocounter{theorem}{1}
\begin{equation}
\label{Enu}
\vecnu_\ji(x)
:=
\A_j((0,0,x)+\veceperp_{i,\alpha_j(x)})
\in T_{\al_j(x)}\N
\end{equation}

\subsection*{The construction of the supports of the wings}
$\phantom{ab}$
\nopagebreak

We intend to construct the initial surface minus the Scherk cores as a graph over an
appropriate perturbation of $\Xhat(\What)$.
The strategy is to first perturb 
$\left.\Xhat\right|_{\Chat}$,
and then use this as boundary data to perturb $\Xhat(\What)$.
$\Xhat_\zzero:\What\to\N$ will be an appropriate parametrization of the perturbed $\Xhat(\What)$ which is defined as
follows.
We first define $\Xhat_\zzero$ on $\Chat$ by requiring (compare with \ref{EEjzero} and recall \ref{DAj} and \ref{PScherk}.ii)
\addtocounter{theorem}{1}
\begin{equation}
\label{EXhat0C}
\Xhat_\zzero\circ\alhat_{j,i}(x)
=
\E_j
\left(\bb_{i,\alpha_j(x)}(\ttt_j^{-1}(x),0)\right),
\quad\text{where}\quad
\ttt_j(x):=\rho_j\, x.
\end{equation}

Given $x\in\Rbb$,
we can determine uniquely a point $\beche_{j,i}(x)\in\W$ in the vicinity of 
$\QQ\circ\alhat_{j,i}(x)$,
and a (small) vector $\left.\cch\right|_{\beche_{j,i}(x)}\in T_{\X\circ\beche_{j,i}(x) }\N$
orthogonal to $\X(\W)$,
such that
\addtocounter{theorem}{1}
\begin{equation}
\label{Ebehat}
\Xhat_\zzero\circ\alhat_{j,i}(x)
=
\exp_{\X\circ\beche_{j,i}(x)}(\left.\cch\right|_{\beche_{j,i}(x) }).
\end{equation}
Note that we have defined a curve $\beche_{j,i}:\Rbb\to\W$,
and a section $\cch$ of the pull-back by $\X\circ\beche_{j,i}$ of the normal bundle of $\X(\W)$,
such that the graph of $\cch$ provides
the point $\Xhat_\zzero\circ\alhat_{j,i}(x)$
above the point $\X\circ\beche_{j,i}(x)$.

We define now $\W_\zzero$ to be the complement of a thin neighborhood of $\C$ in $\W$
such that
\addtocounter{theorem}{1}
\begin{equation}
\label{EW0}
\partial \W_\zzero=\partial W\cup
\left(\cup_{j,i} \,\,\,\beche_{j,i}(\Rbb)\right).
\end{equation}
We can consider then $\cch$ as a normal section on $\partial W_\zzero\setminus\partial\W$.
We extend it to vanish on $\partial W$ (Dirichlet boundary conditions).
Assuming $\cch$ appropriately small,
we can solve the linearized Dirichlet problem and then correct for the nonlinear terms,
so that we obtain an extension of $\cch$ as a normal section to $\W_\zzero$,
uniquely defined by the requirement that the graph of $\cch$ over $\X(\W_\zzero)$ is minimal.
This is of course possible because of the smallness of the perturbations and the flexibility condition
\ref{Tmain}.ii.
$\Xhat_\zzero$ is already defined on $\Chat$ by \ref{EXhat0C}.
We define $\Xhat_\zzero:=\Xhat$ on $\partial\W$.
We extend it then on the whole of $\What$ so that it is a small perturbation of $\Xhat$,
and provides a parametrization of the minimal graph of $\cch$ over $\X(\W_\zzero)$.
(Its precise definition is not needed so we do not provide one).

\addtocounter{equation}{1}
\begin{definition}
\label{Dbnu}
We define smooth maps
$\bbeha_{j,i}:\Rbb\times[0,5\delta_1]\to\What$,
$\bbe_{j,i}:\Rbb\times[0,5\delta_1]\to\Xhat_\zzero(\What)\subset\N$,
and
$\vecmu_{j,i}:\Rbb\times[0,5\delta_1]\to T\N$,
for some small $\delta_1$ to be determined independently of the $m_j$'s,
by requiring the following.
$\bbeha_{j,i}(x,.):[0,5\delta_1]\to\What$ for any fixed $x$
is a unit speed geodesic with respect to the metric induced by $\Xhat_\zzero:\What\to\N$ on $\What$
with initial conditions 
(recall \ref{EXhat0C})
$$
\bbeha_{j,i}(x,0)= \alhat_{j,i}(x) \in \Chat,
$$
and corresponding initial velocity 
equal to the inward unit conormal in the same metric.
$\bbe_{j,i}:=\Xhat_\zzero\circ\bbeha_{j,i}$,
and 
$\vecmu_{j,i}(x,s)$ is the unit normal to the image of $\bbe_{j,i}$ at $\bbe_{j,i}(x,s)$
with orientation chosen so that its coordinates in a geodesic coordinate system around
$\al(x)$ are close to the ones of $\vecnu_\ji(x)$
(recall \ref{Enu}).
\end{definition}

Note that
the $\bbeha_{j,i}$'s provide then covering parametrizations of the tubular neighborhoods of the
components of $\Chat$ in $\What$
and similarly
the $\bbe_{j,i}$'s provide covering parametrizations of the tubular neighborhoods of the
components of $\Xhat_\zzero(\Chat)$ in $\Xhat_\zzero(\What)$.
The $\vecmu_\ji$'s provide the corresponding unit normal fields to the $\bbe_{j,i}$'s.

\subsection*{The construction of the initial surface}
$\phantom{ab}$
\nopagebreak

Recall that parametrizations $\BB_{i,\alpha}$ of the wings of the Scherk surfaces are defined in \ref{PScherk}
as graphs over the parametrizations $\bb_{i,\alpha}$ over the asymptotic half-planes by the functions $\fa$.
We would like to construct the wings of the Scherk surfaces in the initial surface by replacing the parametrizations
$\bb_\ialpha$ of the asymptotic half-planes with the $\bbe_{j,i}$'s we just constructed.
Some difficulties we have to face include 
the scaling by $\rho_j$ we have to introduce,
and the need to transit to different descriptions of the surfaces
at the beginning of the wings (close to the core as in \ref{EEjzero}),
and the end of the wings (by $\Xhat_\zzero$ far away from the core).
There are similarities with the bending of the wings in 
\ref{DZalthph}.iii.

\addtocounter{equation}{1}
\begin{definition}
\label{Dbbeti}
We define 
$\bbeti_{j,i}:\Rbb\times[0,5\delta_1]\to\N$
as follows:
\newline
(i).
On $\Rbb\times[0,1/10]$ we have
$
\Zti_j\circ \BB_{i,\pi/4}=\bbeti_{j,i}\circ\cc_j,
$
where $\cc_j:\Rtwoplus\to\Rtwoplus$ is defined by
$
\cc_j(x,s):=\rho_j\,\,(x,s).
$
\newline
(ii).
On $\Rbb\times[\rho_j,5\delta_1]$
we define 
$\bbeti_{j,i}$ by
\addtocounter{theorem}{1}
\begin{equation}
\label{Ebbeti}
\bbeti_{j,i} (x,s)=\exp_{\bbe_{j,i}(x,s)}
\left(
\ffti_{j,i}(x,s)\,\vecmu_{j,i}(x,s)
\right),
\end{equation}
where $\ffti_{j,i}:\Rtwoplus\to\Rbb$ is defined by
\addtocounter{theorem}{1}
\begin{equation}
\label{Efftipp}
\ffti_{j,i}\circ\cc_j(x,s)=
\rho_j\,
\boldsymbol{f}_{\alpha_j(x)} ( x, s)
\,
\psi[5\delta_1,4\delta_1](\rho_j s).
\end{equation}
(iii).
On $\Rbb\times[\rho_j/10,\rho_j]$
we have to transit from (i) to (ii).
Let $\Phi$ be a local parametrization in the vicinity of 
$\bbe_{j,i}(x,s)$ defined by
\addtocounter{theorem}{1}
\begin{equation}
\label{EPhi}
\Phi(x',s',\chi')=
\exp_{\bbe_{j,i}(x',s')}
(
\chi'\,\vecmu_{j,i}(x',s')
).
\end{equation}
We define then
$\bbeti_{j,i}$ on $\Rbb\times[\rho_j/10,\rho_j]$ 
by (recall \ref{Lpsiab})
\addtocounter{theorem}{1}
\begin{equation}
\label{Ebbeti3}
\Phi^{-1}\circ\bbeti_{j,i}(x,s)=
\psi[\rho_j,0](s)\Phi^{-1}\circ
\Zti_j\circ \BB_{i,\pi/4}\circ\cc_j^{-1}+
\psi[0,\rho_j](s)(x,s,\ffti_{j,i}(x,s)\,)
\end{equation}
\end{definition}

Note that the choice of orientation of
$\vecmu_{j,i}(x,s)$ 
in \ref{Dbnu}
ensures the 
smallness of the third coordinate under $\Phi^{-1}$ of
$\Zti_j\circ \BB_{i,\pi/4}\circ\cc_j^{-1}-(x,s,\ffti_{j,i}(x,s)\,)$
compared to $\ffti_{j,i}(x,s)$.
Now that the wings have been defined, we extend $\Zti_j$ to the ``extended core'' $\Sext_{,j}$
of $\Sc$, where
\addtocounter{theorem}{1}
\begin{equation}
\label{ESext}
\Sext_{,j}:=
\Scortho\setminus\cup_{i=1}^4\BB_{i,\pi/4}(\Rbb\times[5\rho_j^{-1}\delta_1,\infty)),
\end{equation}
by requiring that on $\Sext_{,j}\setminus\Score$
we have
\addtocounter{theorem}{1}
\begin{equation}
\label{EZtiwings}
\Zti_j\circ\BB_\ipi
=
\bbeti_\ji.
\end{equation}

\addtocounter{equation}{1}
\begin{definition}
\label{Dcentral}
We define the central initial surface by
$$
\Mini_\zzero:=
\left(\cup_j\Zti_j(\Sext_{,j})\right)
\cup
\left(\X_\zzero(\What)\setminus \cup_{i,j}
\bbe_{j,i}(\Rbb\times[0,5\delta_1])  \right).
$$
\end{definition}

\section{The family of initial surfaces for the desingularization construction}
\label{Sfamily}
\nopagebreak

\subsection*{The parameters of the family of initial surfaces}
$\phantom{ab}$
\nopagebreak

According to the general methodology we are following,
the parameters of the initial surfaces correspond to the dislocations
we use to create extended substitute kernel.
As we discussed earlier for each hemispherical region of the Scherk surfaces used,
we need to introduce seven dislocations
which are controlled by seven parameters 
which we have denoted by $\sigma$, $\vartheta_1$, $\vartheta_2$, and $\varphi_i$ where $\iequal$.
The construction of our initial surfaces will have therefore $7\sum_{j=1}^k m_j$
parameters.
More precisely we assume given
\addtocounter{theorem}{1}
\begin{equation}
\label{Exibo}
\begin{aligned}
&\xibo:= (\underline{\sigma}_j,\underline{\vartheta}_j,\underline{\varphi}_j)_{j=1}^k,
\quad\text{where}\quad
\\
&\underline{\sigma}_j :=  \{ {\sigma}_\jq \}_{q\in\Zbb},
\quad
\underline{\vartheta}_j := \{ {\vartheta}_\jq \}_{q\in\Zbb},
\quad
\underline{\varphi}_j := \{ {\varphi}_\jq \}_{q\in\Zbb},
\quad\text{where}\quad
\\
&\sigma_{\jq}\in\Rbb,
\quad
\vartheta_{\jq}= (\vartheta_{\jq,1},\vartheta_{\jq,2}) \in\Rtwo,
\quad
\varphi_{\jq}= (\varphi_{\jq,i})_\iequal\in\Rfour.
\end{aligned}
\end{equation}
At the moment we assume that these sequences are such that the constructions that follow are well defined.
In the next section we will specify the range of values allowed for the parameters.
Because of the periodicity involved as for example in \ref{EAlk} and \ref{EEjrho},
we assume appropriate periodicity conditions which reduce the parameter count to $7\sum_{j=1}^k m_j$
as above.
For example we always have to assume $\sigma_{\jq}=\sigma_{j,q+m_j}$,
and for the simplest topology $\varphi_{\jq}=\varphi_{j,q+m_j}$.

In the construction of the initial surfaces that follow
we often have to convert the discrete data provided by the above sequences,
to smooth functions describing attributes which have to vary smoothly along the surfaces.
This motivates the following definition:

\addtocounter{equation}{1}
\begin{definition}
\label{DPsi}
Given a sequence $\underline{v}=\{v_q\}_{q\in\Zbb}$ with values $v_q$ in some vector space $V$,
we define a smooth function $v=\Psi\underline{v}:\Rbb\to V$ by
$v=v_q$ on $[q\pi+1,(q+1)\pi-1]$
and
$v=v_q+(v_{q+1}-v_q)\psi[(q+1)\pi-1,(q+1)\pi+1]$ on $[(q+1)\pi-1,(q+1)\pi+1]$.

We also define a function $\Psi_0\underline{v}:\Rbb\to V$ by $\Psi_0\underline{v}=v_q$ on $(q\pi,(q+1)\pi]$.
\end{definition}

Note that $\Psi_0$ is an operator which converts sequences to step functions with the same values,
and $\Psi$ can be considered as $\Psi_0$ followed by a smoothing.
We define for future reference
for a sequence $\underline{v}$ as above
$\|\underline{v}:\ell^\infty(\Cunder_j)\|:=\max_p|v_p|$,
and for $r\in[1,\infty)$
\addtocounter{theorem}{1}
\begin{equation}
\label{Elrnorm}
\|\underline{v}:\ell^r(\Cunder_j)\|:=
\|\Psi_0\underline{v}:\, L^2([0,m_j\pi],\rho_j^2 \,dx^2)\|=
\left(\frac{\ell_j}{m_j} \sum_{q=1}^{m_j} v_q^r \right)^{1/r}.
\end{equation}

\subsection*{Scaling}
$\phantom{ab}$
\nopagebreak

We assume now a $\xibo$ as above fixed,
and we proceed to modify the construction of $\Mini_\zzero$ carried out in the previous
section, to the construction of $\Mini'=\Mini_\xibo$.
We start by modifying $\CC_j$, which controls the scaling of the Scherk surfaces used, to $\CC'_j$.
We define 
\addtocounter{theorem}{1}
\begin{equation}
\label{ECCp}
\begin{aligned}
&\CC'_j(x_1,x_2,x_3):=
\left(
\rho_j(x_3)\,x_1,\,\rho_j(x_3)\,x_2,\,\ttt'_j(x_3)
\right),
\\
&\text{where }
\quad
\ttt'_j(x_3):=
\int_0^{x_3}\rho_j(x)dx,
\end{aligned}
\end{equation}
where $\rho_j:\Rbb\to\Rbb$ is not a constant anymore 
and will be determined in terms of $\underline{\sigma}_j$.
Because of the periodicity it has to satisfy
\addtocounter{theorem}{1}
\begin{equation}
\label{Erhoperiod}
\rho_j(x+m_j\pi)\equiv\rho_j(x),
\qquad\qquad
\ttt'_j(m_j\pi)=
\int_0^{m_j\pi}\rho_j(t)dt=\ell_j.
\end{equation}

We require that $\rho_j$ is determined by the ODE
\addtocounter{theorem}{1}
\begin{equation}
\label{ErhoODE}
\frac {d\rho_j}{dx}=\sigma_j\rho_j^2,
\quad
\text{ where }
\sigma_j:=\Psi(\underline{\sigma}_j).
\end{equation}
This ODE is motivated by the following heuristic argument:
The amount of kernel created on a half-handle scaled to unit size
is the change $\rho_j((q+1)\pi)-\rho_j(q\pi)$ (which can be approximated with
$\pi\frac {d\rho_j}{dx}$)
divided by the scaling factor
${\rho_j(q\pi)}$.
The amount of kernel needed is proportional to ${\rho_j(q\pi)}$,
and we choose the proportionality constant to be $\sigma_\jq\pi$.
By smoothing and approximating the equation follows.

By rearranging \ref{ErhoODE} and integrating we obtain that 
\addtocounter{theorem}{1}
\begin{equation}
\label{ErhoODEsol}
\rho_j(x)=\left({(\rho}_j(0))^{-1}-\int_0^x\sigma(t)dt\right)^{-1},
\end{equation}
where $\rho_j(0)$ is uniquely determined by the second equation in \ref{Erhoperiod},
while the first amounts then to the condition
\addtocounter{theorem}{1}
\begin{equation}
\label{Erhocondition}
\sum_{q=1}^{m_j}
\sigma_\jq
=0.
\end{equation}
It is possible to check that $\rho_j>0$ and satisfies
\addtocounter{theorem}{1}
\begin{equation}
\label{Erhobound}
\max_x \rho_j(x)
\le
\,
C\,\min_x \rho_j(x)\,
\end{equation}
where $C$
depends only on upper bounds for $\ell_j$ and $\|\underline{\sigma}_j:\ell^\infty(\Cunder_j)\|$
and does not depend on $m_j$.

\subsection*{Transverse unbalancing}
$\phantom{ab}$
\nopagebreak

The creation of transverse unbalancing amounts to perturbing the minimal pieces attached to
the curve of intersection $\Cunder$ so that the opposing conormals are not opposite anymore,
but they rather form a prescribed angle close to $\pi$.
The angle varies along $\Cunder$ and is controlled by the $\underline{\vartheta}_j$'s.
By ignoring the nonlinear terms the perturbed minimal pieces can be considered as graphs over $\X(\W)$
satisfying the linearized equation.
Although we do not need to construct exactly the perturbed minimal pieces yet,
we do need to solve the linearized equation for these graphs in order to determine
the perturbed position of $\Cunder$ which ensures (approximately)
the required transverse unbalancing.
Later on we will construct exactly the graphs which determine the supports of the wings.
Those graphs are further perturbations of the graphs we are implicitly considering now.

\addtocounter{equation}{1}
\begin{definition}
\label{DApj}
Given an appropriately small section $\cch$ of the pull-back by $\left.\X\right|_\C$,
of the normal to $\X(\W)$ bundle,
we define smooth perturbations 
$\al'_j:\Rbb\to\N$
and
$\A'_j:\Rbb^3\to N\Cunder'_j$
of the maps 
$\A_j$ and $\al_j$
defined in \ref{DAj},
where
$\Cunder'_j:=\al'_j(\Rbb)$
and $N\Cunder'_j$ is its normal bundle in $\N$,
as follows.

We define $\al'_j(t):=\Sigma_{12}\cap\Sigma_{23}\cap\Sigma_{13}$,
where 
$\Sigma_{12}$ is the image under the exponential map of a neighborhood of the origin in
$N_{\al_j(x)}\Cunder$
(the fiber of the normal bundle of $\Cunder$ above $\al_j(x)$),
and $\Sigma_{23}$ and $\Sigma_{13}$
are respectively for $i=1,2$
the parallel surfaces to the images under $\X$ of neighborhoods
of $p_i:=\QQ\circ\alhat_\ji(x)=\QQ\circ\alhat_{j,-i}(x)$
(recall \ref{DWhat} and \ref{Evectors}),
and at such a height that they contain
$\exp_{\X(p_i)} (\left.\cch\right|_{p_i} )$.

$\A'_j$ 
is then defined by
$\A'_j(0,0,t):=\al'_j(t)$ and 
the following requirements.
$\A'_j$ is a bundle isometry
and for each $t\in\Rbb$ there is $\alpha''_j(t)$ close to $\alpha_j(t)$
such that 
$\A'_j(R_i(\cos\alpha''_j(t),\sin\alpha''_j(t), t))$
is for $i=\pm1$ (recall \ref{NR}) exactly tangent to $\Sigma_{23}$ 
and for $i=\pm2$ exactly tangent to $\Sigma_{13}$.
Moreover in geodesic coordinates around $\al'_j(t)$ it is a small perturbation of 
$\A_j(R_i(\cos\alpha_j(t),\sin\alpha_j(t), t))$.
\end{definition}

By appealing to the flexibility condition \ref{Tmain}.ii,
we extend now uniquely $\cch$,
which is only defined currently on $\C$,
to a section of the pull-back by $\X$ of the normal bundle of $\X(\W)$ in $\N$
on the whole of $\W$,
by requiring the following:
\newline
(i).
The extended $\chi$ vanishes on $\partial\W$.
\newline
(ii). It is continuous on $\W$.
\newline
(iii). It is smooth on $\W\setminus\C$ where it also satisfies the linearized equation
$\Lop\cch=0$.

Equivalently the pull-back by $\QQ$ of $\cch$ to $\What$,
which we will denote by $\cchhat$,
solves the Dirichlet problem on $\What$ for $\Lop$ and the corresponding Dirichlet boundary data.
$\cchhat$ is then smooth on $\What$,
and we define $\thetahat_\ji:\Rbb\to\Rbb$
by (recall \ref{DAj}.iv and \ref{Enu})
\addtocounter{theorem}{1}
\begin{equation}
\label{Ethetahat}
\veceta_\ji(x)\cdot\cchhat
=
\tan\circ\,\thetahat_\ji(x)
\,\vecnu_\ji(x),
\end{equation}
where the left hand side is the directional derivative of $\cchhat$ along $\veceta_\ji(x)$.
$\thetahat_\ji$ determines then approximately the angle by which the conormal $\X_*(\veceta_\ji(x))$
has to turn when $\Xhat(\What)$ is modified according to $\cchhat$ in $\N$.
We define $\theta_\ji:\Rbb\to\Rbb$ by $\theta_\ji:=\thetahat_\ji\circ \ttt'_j$
(recall \ref{ECCp}).
Recall that our construction of the $\theta_\ji$'s
assumed given $\cch$ on $\C$.
We uniquely determine now $\cch$ by appealing to the unbalancing condition \ref{Tmain}.i
and requiring that for $i=1,2$,
\addtocounter{theorem}{1}
\begin{equation}
\label{Etheta}
\theta_\ji
+
\theta_{j,-i}
=
\Psi\underline{\vartheta}_\ji,
\quad\text{where}\quad
\underline{\vartheta}_\ji:=\{\vartheta_{j,p,i}\}_{p\in\Zbb}.
\end{equation}

We define now an extended tetrad (recall \ref{ET})
\addtocounter{theorem}{1}
\begin{equation}
\label{ETjx}
T_j(x):=\{\vecepji(x),\,\Psi\underline{\varphi}_\ji(x)\}_{i=\pm1,\pm2} ,
\quad\text{where}\quad
\underline{\varphi}_\ji:=\{\varphi_{j,p,i}\}_{p\in\Zbb},
\end{equation}
and $\vecepji(x)$ is 
$(\cos\alpha''_j\circ \ttt'_j(x),\sin\alpha''_j\circ \ttt'_j(x),x)$
rotated by an angle
$\theta_\ji(x)$ around the $x_3$-axis and then reflected by $R_i$.
We define then (recall \ref{EZT}) $\Z'_j:\Sc\to\Rthree$,
$\alpha'_j:\Rbb\to\Rbb$,
$\vartheta'_j:\Rbb\to\Rtwo$,
$\varphi'_j:\Rbb\to\Rfour$,
and
$\phi'_j:\Rbb\to\Rbb$ by 
\addtocounter{theorem}{1}
\begin{equation}
\label{EZp}
\Z'_j(p):=Z_{T_j(x_3(p))}(p),
\qquad
Z_{T_j(x)}=Z_{{\alpha'_j(x),\vartheta'_j(x),\varphi'_j(x),\phi'_j(x)}}.
\end{equation}

\subsection*{The construction of the initial surfaces}
$\phantom{ab}$
\nopagebreak

Now that we have the modifications $\A'_j$, $\CC'_j$,
we define in analogy with \ref{EEjzero},
\addtocounter{theorem}{1}
\begin{equation}
\label{EEjq}
\Zti'_j:=\E'_j\circ\Z'_j,
\qquad
\text{where}
\qquad
\E'_j:=\exp\circ\A'_j\circ\CC'_j.
\end{equation}
We define then the cores of the initial surfaces to be the images $\Zti'_j(\Score)$.
We proceed then to define the supports of the wings in an analogous way
as in the construction of $\Mini_\zzero$ in the previous section.
We first define $\Xhat_\xibo$ on $\Chat$ by requiring
\addtocounter{theorem}{1}
\begin{equation}
\label{EXhatxi}
\Xhat_\xibo\circ\alhat_{j,i}(x)
=
\E'_j
\left(\bb_{i,\alpha'_j(x)}({\ttt'_j}^{-1}(x),0)\right).
\end{equation}

We can define now as in the previous section, a curve $\beche_{j,i}:\Rbb\to\W$
and a section $\cch$ of the pull-back by $\X\circ\beche_{j,i}$ of the normal bundle of $\X(\W)$,
by (recall \ref{Ebehat})
\addtocounter{theorem}{1}
\begin{equation}
\label{Ebehatxi}
\Xhat_\xibo\circ\alhat_{j,i}(x)
=
\exp_{\X\circ\beche_{j,i}(x)}(\left.\cch\right|_{\beche_{j,i}(x) }).
\end{equation}
Both $\beche_\ji$ and $\cch$ depend implicitly on $\xibo$ and although
are denoted with the same symbols as the corresponding objects in \ref{Ebehat},
they are actually modifications of them.

As in the previous section we proceed to
define $\W_\xibo$ as the complement of a thin neighborhood of $\C$ in $\W$
such that
\addtocounter{theorem}{1}
\begin{equation}
\label{EWxi}
\partial \W_\xibo=\partial W\cup
\left(\cup_{j,i} \,\,\,\beche_{j,i}(\Rbb)\right).
\end{equation}
We can consider then $\cch$ as a normal section on $\partial W_\xibo\setminus\partial\W$
and extend it as in the previous section to the whole of $\W_\xibold$
so that its graph is minimal.
We extend then 
$\Xhat_\xibold$ to the whole of $\What$ so that it is a small perturbation of $\Xhat$,
and provides a parametrization of the minimal graph of $\cch$ over $\X(\W_\xibold)$.
We define now maps which depend implicitly on $\xibo$ by essentially repeating \ref{Dbnu}:

\addtocounter{equation}{1}
\begin{definition}
\label{Dbnuxi}
We define smooth maps
$\bbeha_{j,i}:\Rbb\times[0,5\delta_1]\to\What$,
$\bbe_{j,i}:\Rbb\times[0,5\delta_1]\to\Xhat_\xibo(\What)\subset\N$,
and
$\vecmu_{j,i}:\Rbb\times[0,5\delta_1]\to T\N$,
by requiring the following.
$\bbeha_{j,i}(x,.):[0,5\delta_1]\to\What$ for any fixed $x$
is a unit speed geodesic with respect to the metric induced by $\Xhat_\xibo:\What\to\N$ on $\What$
with initial conditions 
(recall \ref{EXhatxi})
$
\bbeha_{j,i}(x,0)= \alhat_{j,i}(x) \in \Chat,
$
and corresponding initial velocity 
equal to the inward unit conormal in the same metric.
$\bbe_{j,i}:=\Xhat_\xibo\circ\bbeha_{j,i}$,
and 
$\vecmu_{j,i}(x,s)$ is the unit normal to the image of $\bbe_{j,i}$ at $\bbe_{j,i}(x,s)$
with orientation chosen so that its coordinates in a geodesic coordinate system around
$\al(x)$ are close to the ones of $\vecnu_\ji(x)$
(recall \ref{Enu}).
\end{definition}

In analogy then with \ref{Dbbeti} we have the following:

\addtocounter{equation}{1}
\begin{definition}
\label{Dbbetixi}
We define 
$\bbeti_{j,i}:\Rbb\times[0,5\delta_1]\to\N$
as follows:
\newline
(i).
On $\Rbb\times[0,1/10]$ we have
$
\Zti'_j\circ \BB_{i,\pi/4}=\bbeti_{j,i}\circ\cc'_j,
$
where $\cc'_j:\Rtwoplus\to\Rtwoplus$ is defined by
$
\cc'_j(x,s):=(\ttt'_j(x),\,\rho_j(x)\,s).
$
\newline
(ii).
On $\cc'_j(\Rbb\times[1,\infty)\,)\cap\Rbb\times[0,5\delta_1]$
we define 
$\bbeti_{j,i}$ by
\addtocounter{theorem}{1}
\begin{equation}
\label{Ebbetixi}
\bbeti_{j,i} (x,s)=\exp_{\bbe_{j,i}(x,s)}
\left(
\ffti_{j,i}(x,s)\,\vecmu_{j,i}(x,s)
\right),
\end{equation}
where $\ffti_{j,i}:\Rtwoplus\to\Rbb$ is defined by
\addtocounter{theorem}{1}
\begin{equation}
\label{Efftippxi}
\ffti_{j,i}\circ\cc'_j(x,s)=
\rho_j(x)\,
\boldsymbol{f}_{\alpha'_j(x)} ( x, s)
\,
\psi[5\delta_1,4\delta_1](\rho_j(x)\, s).
\end{equation}
(iii).
Let $\Phi$ be a local parametrization in the vicinity of 
$\bbe_{j,i}(x,s)$ defined by
\addtocounter{theorem}{1}
\begin{equation}
\label{EPhixi}
\Phi(x',s',\chi')=
\exp_{\bbe_{j,i}(x',s')}
(
\chi'\,\vecmu_{j,i}(x',s')
).
\end{equation}
We define then
$\bbeti_{j,i}$ on
$\cc'_j(\Rbb\times[1/10,1]\,)$ 
by (recall \ref{Lpsiab})
\addtocounter{theorem}{1}
\begin{equation}
\label{Ebbeti3xi}
\Phi^{-1}\circ\bbeti_{j,i}(x,s)=
\psi[\rho_j,0](s)\Phi^{-1}\circ
\Zti'_j\circ \BB_{i,\pi/4}\circ{\cc'_j}^{-1}+
\psi[0,\rho_j](s)(x,s,\ffti_{j,i}(x,s)\,)
\end{equation}
\end{definition}

We modify \ref{ESext} to
\addtocounter{theorem}{1}
\begin{equation}
\label{ESextxi}
\Sext_{,j}:=
\Scortho\setminus\cup_{i=1}^4\BB_{i,\pi/4}({\cc'_j}^{-1}(\Rbb\times[0,5\delta_1]).
\end{equation}
As in the previous section we extend $\Zti_j$ to the extended core $\Sext_{,j}$,
by requiring that on $\Sext_{,j}\setminus\Score$
we have
\addtocounter{theorem}{1}
\begin{equation}
\label{EZtiwingsxi}
\Zti'_j\circ\BB_\ipi
=
\bbeti_\ji.
\end{equation}

\addtocounter{equation}{1}
\begin{definition}
\label{DMxi}
For $\xibo$ as in \ref{Exibo} satisfying \ref{Erhocondition}
we define the corresponding initial surface by
$$
\Mini_\xibo:=
\left(\cup_j\Zti'_j(\Sext_{,j})\right)
\cup
\left(\X_\xibo(\What)\setminus \cup_{i,j}
\bbe_{j,i}(\Rbb\times[0,5\delta_1])  \right).
$$
\end{definition}

\section{Main estimates and outline of the proof}
\label{Smain}
\nopagebreak

\subsection*{The range of the parameters}
$\phantom{ab}$
\nopagebreak

In this section we discuss some estimates and give a rough outline of the proof.
We start by discussing the range of the parameters.
Recall \ref{Exibo}.
We define $\Xi$ to be the vector space of sequences as $\xibo$ which are subject to the
appropriate periodicity conditions as discussed in the previous section.
We also define the subspace $\Xi_0$ to be those sequences in $\Xi$ which also satisfy
\ref{Erhocondition}.
We have then
\addtocounter{theorem}{1}
\begin{equation}
\label{EdimXi}
\dim\Xi=
7\sum_{j=1}^k m_j,
\qquad\qquad
\dim\Xi_0=
7\sum_{j=1}^k m_j-k.
\end{equation}

Recall \ref{Cmj}.
We define
$\tau:=1/m_1$ and we have then
\addtocounter{theorem}{1}
\begin{equation}
\label{Etaumj}
\cunder^{-1}_0\tau^{-1}\le m_j \le \cunder_0\tau^{-1}.
\end{equation}

We will need a kind of discrete derivative for a sequence $\underline{v}$ as the one in \ref{Elrnorm},
and therefore we define 
\addtocounter{theorem}{1}
\begin{equation}
\label{Edd}
\dd\underline{v}:=\{\tau^{-1}(v_{q+1}-v_q)\}_{q\in\Zbb}.
\end{equation}
It would be more accurate to use the factor $m_j$ instead of $\tau^{-1}$,
but because of \ref{Etaumj} that would modify $\dd$ only by a factor controlled by $\cunder_0$
and therefore would not make a difference in our presentation.
We also define for $k\in\Nbb_0$ and $r\in[1,\infty]$ (recall \ref{Elrnorm})
\addtocounter{theorem}{1}
\begin{equation}
\label{Elkrnorm}
\|\underline{v}:\ell_k^r(\Cunder_j)\|:=
\sum_{k'=0}^k
\|\dd^{k'}\underline{v}:\ell^r(\Cunder_j)\|.
\end{equation}

We define now the range of the parameters $\xibold$ by requiring 
\addtocounter{theorem}{1}
\begin{equation}
\label{Erange}
\|\underline{\sigma}_j:\ell_1^2(\Cunder_j)\|\le \cunder_1,
\qquad
\|\underline{\vartheta}_j:\ell_1^2(\Cunder_j)\|\le \cunder_2,
\qquad
\|\underline{\varphi}_j:\ell_1^2(\Cunder_j)\|\le \cunder_2,
\end{equation}
where $\cunder_1$ and $\cunder_2$ are chosen later independently of $\tau$
but large enough depending on the given system of minimal surfaces and $\cunder_0$.

\subsection*{The mean curvature of the initial surfaces}
$\phantom{ab}$
\nopagebreak

Estimating the mean curvature is the final product of a long process
where one has to estimate many error terms.
In particular $\cch$ has to be estimated carefully.
We only outline here the final result.
Note that by the definitions of the initial surfaces \ref{Dcentral}
and \ref{DMxi} the mean curvature is supported on the extended cores.
In order to discuss the estimates it is helpful to define components of the mean
curvature such that
each of them has support contained in the support of the mean curvature and
\addtocounter{theorem}{1}
\begin{equation}
\label{EHparts}
H = \Hgluing+\Herror+\Hprescribed,
\end{equation}
where each of the components (that is the summands on the right)
satisfies different estimates:
In some abstract sense $\Hgluing$ is created by the gluing construction
and satisfies estimates independent of the size of $\cunder_1$ and $\cunder_2$,
and $\Herror+\Hprescribed$ is created by the dislocations controlled by $\xibo$
and therefore their estimates depend on $\cunder_1$ and $\cunder_2$.

Before we proceed we remark that the mean curvature and its components above are sections of the normal bundle
of $\Mini_\xibo$ but when we pull them back by 
$\Zti'_j$
to $\Sext_{,j}$, 
we can identify them with functions.
Since these functions are supported away from the boundary of $\Sext_{,j}$
we can extend them smoothly to the whole of $\Sc$ by having them vanish on $\Sc\setminus\Sext_{,j}$:

\addtocounter{equation}{1}
\begin{convention}
\label{CHfunction}
From now on we will consider the pullback by $\Zti'_j$ of the mean curvature and its components
as smooth functions on $\Sc$ as described above.
\end{convention}

We start by describing $\Hprescribed$.
For this we define first a ``model'' or ``tangent'' embedding 
\addtocounter{theorem}{1}
\begin{equation}
\label{EWjq}
\WW_{j,q}:\Sc\cap\{x_3\in[(q-1)\pi,(q+2)\pi]\to\Rthree,
\end{equation}
which depends only on the parameters
$\widetilde{\sigma}_\jq:= \rho_j(q\pi)\, \sigma_{j,q},\, \vartheta_{j,q,i'},\, \varphi_{j,q,i}$,
where $i'=1,2$, $\iequal$.
When $\sigma_{j,q}=0$
we define
$\WW_{j,q}$ to be simply $Z_\althph$ as in \ref{DZalthph}
with $\vartheta=\{\vartheta_{j,q,i'}\}_{i'=1,2}$
and $\varphi=\{\varphi_{j,q,i}\}_\iequal$.
When $\sigma_{j,q}\ne0$ we modify $Z_\althph$ so that it 
remains unchanged on $\{x_3=q\pi\}$, depends smoothly on the parameters,
and has the scale changing appropriately along $x_3$ at a rate of $\widetilde{\sigma}_\jq$.

For $i=0,\pm1,\pm2,3,4$ we define functions $w_{q,i} : \Sc\cap\{x_3\in[(q-1)\pi,(q+2)\pi]\} \to\Rbb$ 
as the linearizations of the mean curvature of $\WW_{j,q}$ with respect to each of the parameters,
that is we have
\addtocounter{theorem}{1}
\begin{equation}
\label{Ewjq}
H\circ\WW_{j,q}
=
\widetilde{\sigma}_{j,q}w_{q,0}
+\vartheta_{j,q,1}w_{q,3}
+\vartheta_{j,q,2}w_{q,4}
+\sum_\iequal \varphi_{j,q,i} w_{q,i}
+O,
\end{equation}
where $O$ is bounded by a constant times the squares of the parameters.
Note that $w_{q,i}$ does not depend on $j$ and its dependence on $q$ is only because 
its domain changes by a translation.
In analogy with \ref{DPsi} we have now the following:

\addtocounter{equation}{1}
\begin{definition}
\label{DPsif}
Given a sequence of smooth functions
$$
\underline{f}=\{f_q\}_{q\in\Zbb},
\quad\text{where}\quad
f_q:\Sc\cap\{x_3\in[(q-1)\pi,(q+2)\pi]\} \to\Rbb,
$$
we define a smooth function $f=\Psi\underline{f}:\Sc \to \Rbb$ by
$f=f_q$ on $\Sc\cap\{x_3\in[q\pi+1,(q+1)\pi-1]\}$,
and
$f=f_q+(f_{q+1}-f_q)\psi[(q+1)\pi-1,(q+1)\pi+1]\circ x_3$ on
$\Sc\cap\{x_3\in[(q+1)\pi-1,(q+1)\pi+1]\}$.
\end{definition}

\addtocounter{equation}{1}
\begin{definition}
\label{DTheta}
Given 
$$
\llambda=
\left(
\{\lambda_{\jq,0}\}_{q\in\Zbb},
\{(\lambda_{\jq,i})_{i=3,4}\}_{q\in\Zbb},
\{(\lambda_{\jq,i})_{\iequal}\}_{q\in\Zbb}
\right)_{j=1}^k
\in\Xi
$$
we define $\Theta(\llambda)$ to be  a normal section on $\Mini_0$ (recall \ref{CHfunction}),
supported on its extended cores,
where it is determined by
$\Theta(\llambda) \circ \Zti'_j =\Psi\underline{f} \vecnu$
on $\Sext_{,j}$, where $\Psi$ is as in \ref{DPsif},
$\underline{f}=\{f_q\}_{q\in\Zbb}$,
and
$f_q:\Sc\cap\{x_3\in[(q-1)\pi,(q+2)\pi]\} \to\Rbb$
is defined by 
$$
f_q=
\left(
\lambda_{\jq,0}w_{q,0}
+
\frac1 {\rho_j(q\pi)}
\sum_{\iequal,3,4} \lambda_{j,q,i} w_{q,i}
\right)
{\widetilde{\psi}_j},
$$
where
$\widetilde{\psi}_j:\Sc\to[0,1]$ is a cut-off function defined by
$\widetilde{\psi}_j\equiv1$
on 
$\Score$,
and for 
$(x,s)\in\Rtwoplus$,
$\iequal$,
$$
\widetilde{\psi}_j \circ  \BB_{i,\pi/4} (x,s)=
\psi[5\delta_1,4\delta_1](\rho_j(x)\, s).
$$
\end{definition}

Note that $\Theta(\llambda)$ depends on $\xibo$ because it is defined on $\Mini_\xibo$
and also $f_q$ above depends on $\xibo$ through $\rho_j$.
We define now
\addtocounter{theorem}{1}
\begin{equation}
\label{EHprescribed}
\Hprescribed:=\Theta(\xibo).
\end{equation}

In order to describe now the estimates for $\Hgluing$ and $\Herror$ we need to define norms
for functions on $\Sc$.
Given a function $f$ on $\Sc$ we first define in analogy with \ref{Edd}
another function $\dd f$ on $\Sc$ by
\addtocounter{theorem}{1}
\begin{equation}
\label{Eddf}
\dd f:=\tau^{-1}(f\circ\tran-f),
\end{equation}
where $\tran:\Sc\to\Sc$ is the restriction to $\Sc$ of the translation in $\Rthree$ by $(0,0,2\pi)$.
In analogy with \ref{Elrnorm}
we define 
\addtocounter{theorem}{1}
\begin{equation}
\label{Elrnormf}
\begin{aligned}
\|\phi:\ell^r(C^{n,\beta},h,f)\|:=
\|\underline{v}:\ell^r(\Cunder_j)\|,
\qquad\text{where}\quad
\\
v_q:=\|f:C^{n,\beta}(\Sc\cap \{x_3\in[(q-1)\pi,(q+2)\pi]\},h,f)\|,
\end{aligned}
\end{equation}
where $h$ is a metric on $\Sc$ and $f$ a weight function as in 
\ref{E:weightedHolder}.
In analogy with \ref{Elkrnorm} we also define
\addtocounter{theorem}{1}
\begin{equation}
\label{Elkrnormf}
\|\phi:\ell^r_k(C^{n,\beta},h,f)\|:=
\sum_{k'=0}^k
\|\dd^{k'}
\phi:\ell^r(C^{n,\beta},h,f)\|.
\end{equation}

The estimates we have then for $\Hgluing$ and $\Herror$ are as follows.
\addtocounter{theorem}{1}
\begin{equation}
\label{EHestimates}
\begin{aligned}
\|\Hgluing
\circ \Zti'_j
:
\ell^2_1(C^{0,\beta},g,e^{-\gamma s})
\|
&\le
C,
\\
\|\Herror
\circ \Zti'_j
:
\ell^2(C^{0,\beta},g,e^{-\gamma s})
\|
&\le
C(\cunder_1,\cunder_2) \,\tau,
\end{aligned}
\end{equation}
where $\beta,\gamma\in (0,1)$ are fixed constants
and  the constants $C$ and $C(\cunder_1,\cunder_2)$
depend on the given system of minimal surfaces and $\cunder_0$,
but they do not depend on $\tau$ and (equivalently) the $m_j$'s.
The first constant does not depend on $\cunder_1$ and $\cunder_2$,
but the second does.

Regarding these estimates we have the following remarks:
$g$ is the induced metric on $\Sc$.
We could be using instead the metric induced by $\Zti'_j$
rescaled so that the Scherk handles have unit size.
Such an estimate would be equivalent to the one we have up to uniform constants.
The weight function 
$e^{-\gamma s}$ is defined on the wings where $s$ is the second coordinate of $\Rtwoplus$ as usual,
and extended to the rest of $\Sc$ (the core without a small margin)
to be $\equiv1$.
Note that it would be appropriate to consider the mean curvature of the rescaled extended cores
so that the Scherk handles are of unit size.
This amounts to a reduction of the mean curvature by factors of the size of $\rho_j$,
which would mean an extra factor of $\tau$ on the bounds.
This is consistent with the estimates in \cite{kapouleas:finite}.
Finally note that the estimate for $\Herror$ makes it negligible compared with $\Hgluing$
at the $\ell^\infty(C^{0,\beta})$ level.
On the other hand although both discrete derivatives $\dd \Hgluing$ and $\dd\Herror$ are of order $\tau$,
$\dd\Herror$ can be much larger than $\dd \Hgluing$ if $\cunder_1$ and $\cunder_2$ are chosen large
because the second constant depends on $\cunder_1$ and $\cunder_2$ while the first does not.

\subsection*{The linearized equation}
$\phantom{ab}$
\nopagebreak

We need to solve now the linearized equation on $\Mini_\xibo$
\addtocounter{theorem}{1}
\begin{equation}
\label{Elinear}
\Lop u= H + \Theta(\llambda),
\end{equation}
where the $w_{q,i}$'s in the term $\Theta(\llambda)$
(recall \ref{DTheta})
play the role of a basis of the extended substitute kernel.
We are solving modulo the image of $\Theta$ and this way we can obtain the required estimates on $u$.
An important special feature of this construction is that estimates on the rate of change of 
the components of $\llambda$ as the handles vary along $\Cunder$ are crucial for closing the argument.
The easier part of the strategy for solving \ref{Elinear}
is to use the $\Zti'_j$'s to transplant the equation to the 
$\Sext_{,j}$'s
and solve the equation there with Dirichlet boundary conditions.
An iteration then in the usual fashion (see \cite{kapouleas:finite} for a similar step)
provides a global solution.

Our strategy is now to solve the equation for each component separately.
By \ref{EHprescribed}
$\Hprescribed$ is already in the image of $\Theta$.
To solve for $\Hgluing$ we follow the following strategy:
We approximate the problem semi-locally,
that is on the domain of 
$\WW_{j,q}$,
using the operator induced by a rescaling of $\WW_{j,q}$.
Because of the available estimate \ref{EHestimates}
on $\dd\Hgluing$,
we can slightly modify $\Hgluing$ and extend to the whole of $\Sc$
so we have periodicity with respect to the translation by $(0,0,2\pi)$.
The extended substitute kernel has then dimension $7$ and we can solve the equation
in a way similar to the linearized equation on the extended cores in \cite{kapouleas:finite}.
We can patch together the semi-local solutions as in \ref{DPsif}.
This provides $u'$ and $\llambda'$ and an error term $E$ 
such that
\addtocounter{theorem}{1}
\begin{equation}
\label{Euone}
\Lop u'= \Hgluing + E + \Theta(\llambda'),
\end{equation}
where
$
\llambda'=
\left(
\{\lambda'_{\jq,0}\}_{q\in\Zbb},
\{(\lambda'_{\jq,i})_{i=3,4}\}_{q\in\Zbb},
\{(\lambda'_{\jq,i})_{\iequal}\}_{q\in\Zbb}
\right)_{j=1}^k
\in\Xi
$
and the following estimates hold:
\addtocounter{theorem}{1}
\begin{equation}
\label{Eupestimates}
\begin{aligned}
\|u'
\circ \Zti'_j
:
\ell^2_1(C^{2,\beta},g,e^{-\gamma s})
\|
&\le
C(\cunder_1)\,\tau^2,
\\
\|E
\circ \Zti'_j
:
\ell^2(C^{2,\beta},g,e^{-\gamma s})
\|
&\le
C(\cunder_1,\cunder_2) \,\tau,
\\
\|\{\lambda'_{\jq,0}\}_{q\in\Zbb}:\ell^2_1(\Cunder_j)\|
&\le C,
\\
\|\{\lambda'_{\jq,i}\}_{q\in\Zbb}:\ell^2_1(\Cunder_j)\|
&\le C(\cunder_1)\,\tau
\quad
(\iequal,3,4).
\end{aligned}
\end{equation}
Note that the $\tau^2$ factor in the first estimate reduces to $\tau$
if we interpret the estimate in the natural scale of the extended Scherk cores.
Similarly the $\tau$ factor in the second estimate would improve to $\tau^2$.

The next step is
to solve the equation
\addtocounter{theorem}{1}
\begin{equation}
\label{Eutwo}
\Lop u''= \Herror - E + \Theta(\llambda''),
\end{equation}
where
$
\llambda''=
\left(
\{\lambda''_{\jq,0}\}_{q\in\Zbb},
\{(\lambda''_{\jq,i})_{i=3,4}\}_{q\in\Zbb},
\{(\lambda''_{\jq,i})_{\iequal}\}_{q\in\Zbb}
\right)_{j=1}^k
\in\Xi
$
and $E$ is as in \ref{Euone}.
To solve this equation we need to understand the small eigenvalues globally on the extended standard regions
in the $h$ metric,
and also the functions $v$ we use according to the methodology we follow
to ensure exponential decay of the solutions.
The $h$ metric is defined as usual by a conformal change of the induced metric and approximates the spherical metric
pulled back by the Gauss map.
A crucial lemma asserts that a linear combination $u$ of eigenfunctions with small enough eigenvalues 
can be approximated in the $L^2(h)$ metric on each hemispherical region 
$
S_q=\Sa\cap\{x_3\in[q\pi,(q+1)\pi]\}
$
by Jacobi fields corresponding to translations by $\vece_q$
and where we have
\addtocounter{theorem}{1}
\begin{equation}
\label{Eddkernel}
\|\dd\{\vece_q\}_{q\in\Zbb}:\ell^2(\Cunder_j)\|\le
\varepsilon \|\{\vece_q\}_{q\in\Zbb}:\ell^2(\Cunder_j)\|
\le 5\varepsilon (\ell_j/m_j)^{1/2}
\|u:L^2(h)\|
.
\end{equation}
A similar lemma applies to the functions $v$ needed for the decay.
Using these lemmas we conclude the following estimates.
\addtocounter{theorem}{1}
\begin{equation}
\label{Euppestimates}
\begin{aligned}
\|u''
\circ \Zti'_j
:
\ell^2(C^{2,\beta},g,e^{-\varepsilon s})
\|
&\le
C(\cunder_1,\cunder_2,\varepsilon)\,\tau^3,
\\
\|\{\lambda''_{\jq,0}\}_{q\in\Zbb}:\ell^2_1(\Cunder_j)\|
&\le C_1(\cunder_1,\cunder_2)\varepsilon,
\\
\|\{\lambda''_{\jq,i}\}_{q\in\Zbb}:\ell^2_1(\Cunder_j)\|
&\le C_1(\cunder_1,\cunder_2)\varepsilon\,\tau
\quad
(\iequal,3,4).
\end{aligned}
\end{equation}

By \ref{EHparts}, \ref{EHprescribed}, \ref{Euone}, and \ref{Eutwo},
we conclude then
\addtocounter{theorem}{1}
\begin{equation}
\label{Eu}
\Lop u + \Theta(\xibo-\llambda)=H,
\quad\text{where} \quad
u:=u'+u'',
\quad
\llambda:=\llambda'+\llambda''\in\Xi.
\end{equation}
By choosing $\varepsilon=1/C_1(\cunder_1,\cunder_2)$ in \ref{Euppestimates}
and using also \ref{Eupestimates} we conclude
\addtocounter{theorem}{1}
\begin{equation}
\label{Euestimates}
\begin{aligned}
\|u
\circ \Zti'_j
:
\ell^2_1(C^{2,\beta},g,e^{-s/C(\cunder_1,\cunder_2)})
\|
&\le
C(\cunder_1\cunder_2)\,\tau^2,
\\
\|\{\lambda_{\jq,0}\}_{q\in\Zbb}:\ell^2_1(\Cunder_j)\|
&\le C,
\\
\|\{\lambda_{\jq,i}\}_{q\in\Zbb}:\ell^2_1(\Cunder_j)\|
&\le C\,\tau
\quad
(\iequal,3,4).
\end{aligned}
\end{equation}

\subsection*{Closing the argument}
$\phantom{ab}$
\nopagebreak

The estimate we have now for $u$ in \ref{Euestimates} is good enough to ensure
that the quadratic terms satisfy improved estimates by a factor of $\tau/C(\cunder_1,\cunder_2)$.
By choosing now $\cunder_1,\cunder_2$ large enough we can ensure the range of $\xibo$ implied by \ref{Erange}
includes the projections to $\Xi_0$
of the $\llambda$'s which satisfy the estimates in
\ref{Euestimates}.
This allows us to use a fixed point theorem as usual (see for example \cite{kapouleas:finite})
to conclude that for one value of $\xibo\in\Xi_0$ satisfying \ref{Erange}
there is a perturbation $M$ of $\Mini_\xibo$
which satisfies
\addtocounter{theorem}{1}
\begin{equation}
\label{EHfinal}
H=\sum_j \mu_j \widetilde{w}_j,
\end{equation}
where $\widetilde{w}_j:=\Theta(\llambda)$
where $\llambda\in\Xi$ is as in \ref{DTheta} with all the entries vanishing except for $\lambda_{\jq,0}\equiv1$.

It remains to prove that we can arrange for the coefficients $\mu_j$ to vanish.
For this we use that the construction so far has $k$
(recall \ref{ECj})
free continuous parameters corresponding to the choice of 
$\al_j(0)$ in \ref{DAj}.v.
By using the information we already have we can prove that $M$ has
no zero eigenvalues except perhaps the ones corresponding to the eigenfunctions
corresponding to the $\widetilde{w}_j$'s.
By continuously varying the points $\al_j(0)$ we can use the implicit function theorem to obtain a smooth
family of surfaces satisfying \ref{EHfinal}.
By ensuring uniform estimates based on the constructions we already have,
we can prove
periodicity of the family by varying each $\al_j(0)$ along $\Cunder_j(0)$
by a distance of order $\tau$ corresponding to a handle of $\Mini_\xibo$.
Equivalently we can consider the domain of the parameters of the family to be topologically the product of $k$
circles, one circle for each component of $\Cunder$.
By considering then a surface of maximal area in the family,
and using the vanishing of its variations within the family,
we conclude the orthogonality of the mean curvature to certain functions close to $\widetilde{w}_j$.
This implies the vanishing of the $\mu_j$'s in \ref{EHfinal},
and hence the minimality of the corresponding surface.

\bibliographystyle{amsplain}
\bibliography{paper}
\end{document}